\documentclass[bj,preprint]{imsart}

\usepackage{amsthm,amsmath,amssymb,natbib}
\RequirePackage[colorlinks,citecolor=blue,urlcolor=blue]{hyperref}

\newtheorem{theorem}{Theorem}

\newtheorem{proposition}[theorem]{Proposition}

\newtheorem{corollary}[theorem]{Corollary}

\newtheorem{remark}[theorem]{Remark}
\let\oldremark\remark
\renewcommand{\remark}{\oldremark\normalfont}
\newtheorem{example}{Example}
\let\oldexample\example
\renewcommand{\example}{\oldexample\normalfont}

\newcommand{\C}{\mathbb{C}}                       % complex numbers
\newcommand{\E}{\mathop{\mathbb{E}}}              % Expectation operator
\renewcommand{\H}{\mathcal{H}}                    % Hilbert space
\renewcommand{\hat}{\widehat}                     % make hat wide by default
\newcommand{\inv}{^{-1}}                          % inverse operator
\renewcommand{\L}{\mathcal{L}}                    % Lebesgue space
\newcommand{\MSE}{\operatorname{MSE}}             % mean squared error
\newcommand{\N}{\mathbb{N}}                       % naturals
\renewcommand{\P}{\mathcal{P}}                    % family of distributions
\newcommand{\R}{\mathbb{R}}                       % reals
\newcommand{\sminus}{\backslash}                  % set difference
\newcommand{\Var}{\mathbb{V}}                     % variance
\newcommand{\W}{\mathcal{W}}                      % General Sobolev spaces
\newcommand{\X}{\mathcal{X}}                      % base space
\newcommand{\Z}{\mathcal{Z}}                      % index space
\renewcommand{\tilde}{\widetilde}                 % make tilde wide by default
\newcommand{\B}{\mathbb{B}}                       % bias
\DeclareMathOperator{\sinc}{sinc}                 % Sinc function

\begin{document}

\begin{frontmatter}

% "Title of the Paper"
\title{Minimax Estimation of Quadratic Fourier Functionals}

\runtitle{Estimating Quadratic Fourier Functionals}

\begin{aug}
% indicate corresponding author with \corref{}
% \author{\fnms{John} \snm{Smith}\thanksref{a}\corref{}\ead[label=e1]{smith@foo.com}\ead[label=e2,url]{www.foo.com}}
% \address[a]{\printead{e1};\printead{e2}}

\author{\fnms{Shashank} \snm{Singh}\thanksref{CMU}\corref{}\ead[label=e1]{sss1@cs.cmu.edu}}
\author{\fnms{Bharath K.} \snm{Sriperumbudur}\thanksref{PSU}\ead[label=e2]{bks18@psu.edu}}
\and
\author{\fnms{Barnab\'as} \snm{P\'oczos}\thanksref{CMU}\ead[label=e3]{bapoczos@cs.cmu.edu}}

\address[CMU]{Carnegie Mellon University, Pittsburgh, PA 15213, USA
\printead{e1,e3}}

\address[PSU]{Pennsylvania State University, University Park, PA 16802, USA
\printead{e2}}

\runauthor{Singh, Sriperumbudur, and P\'oczos}

\affiliation{Carnegie Mellon University and Pennsylvania State University}

\end{aug}

\begin{abstract}
We study estimation of (semi-)inner products between two nonparametric probability distributions, given IID samples from each distribution.
These products include relatively well-studied classical $\L^2$ and Sobolev inner products, as well as those induced by translation-invariant reproducing kernels, for which we believe our results are the first.
We first propose estimators for these quantities, and the induced (semi)norms and (pseudo)metrics. We then prove non-asymptotic upper bounds on their mean squared error, in terms of weights both of the inner product and of the two distributions, in the Fourier basis.
Finally, we prove minimax lower bounds that imply rate-optimality of the proposed estimators over Fourier ellipsoids.
\end{abstract}

%\begin{keyword}
%\kwd{}
%\kwd{}
%\end{keyword}

% history:
% \received{\smonth{1} \syear{0000}}

%\tableofcontents

\end{frontmatter}

\section{Introduction}
\label{sec:intro}

Let $\X$ be a compact subset of $\R^D$ endowed with the Borel $\sigma$-algebra and let 
%Fix a compact sample space $\X \subseteq \R^D$. Let 
$\P$ denote the family of all Borel probability measures on $\X$. For each $P \in \P$, let $\phi_P : \R^D \to \C$ denote the characteristic function of $P$ given by
% \Bk{I think $\phi_P$ is bit of non-standard notation for characteristic function. May be we should we $\phi_\P$ and accordingly use the different symbol for the Fourier basis.} \CC{I used $\phi_P$ as an alternative to $\hat P$, which is common in functional analysis for the Fourier transform, but is used here as the estimate of $\phi_P$. How about $\phi_P$ for the characteristic function, $\hat \phi_P$ for its estimate, and, use $g_z$ for the Fourier basis?}\Bk{Sounds good if it is not too much of work changing stuff.}
\begin{equation}
\phi_P(z) = \E_{X \sim P} \left[ \overline{\psi_z(X)} \right]
  \quad \text{ for all } z \in \R^D,
  \quad \text{ where } \quad
  \psi_z(x) = \exp \left( i \langle z, x \rangle \right)
  \label{def:fourier_coefficients}
\end{equation}
denotes the $i^{th}$ Fourier basis element, in which $\langle \cdot, \cdot \rangle$ denotes the Euclidean inner product on $\R^D$.
% For any $z \in \mathbb{Z}^D$, let $\psi_z : \X^D \to \C$, defined by
% \[\psi_z(x) = \exp \left( i \langle z, x \rangle \right)\]
% (where $\langle \cdot, \cdot \rangle$ denotes the usual Euclidean inner product on $\R^D$) for all $x \in \X$, denote the usual $z^{th}$ Fourier basis element. For any probability measure $P \in \P$ and $z \in \mathbb{Z}^D$, let
% \[\phi_P(z) := \E_{X \sim P} \left[ \psi_z(X) \right]\]
% denote the $z^{th}$ Fourier coefficient of $P$ (which is always well-defined, since $\|\psi_z\|_\infty \leq 1$).

For any family $a = \{a_z\}_{z \in \Z} \subseteq \R$ of real-valued coefficients indexed by a countable set $\Z$, define a set of probability measures
\[\H_a := \left\{ P \in \P : \sum_{z \in \Z} \frac{\left| \phi_P(z) \right|^2}{a_z^2} < \infty \right\}.\]
Now fix two unknown probability measures $P, Q \in \H_a$. We study estimation of the semi-inner product\footnote{For a complex number $\xi = a + bi \in \C$, $\overline{\xi} = a - bi \in \C$ denotes the complex conjugate of $\xi$. A semi-inner product has all properties of an inner product, except that $\langle P, P \rangle = 0$ does not imply $P = 0$.}
\begin{equation}
\langle P, Q \rangle_a
  = \sum_{z \in \Z} \frac{\phi_P(z) \overline{\phi_Q(z)}}{a_z^2},
\label{eq:estimand}
\end{equation}
as well as the squared seminorm $\|P\|_a^2 := \langle P, P \rangle_a$ and squared pseudometric $\|P - Q\|_a^2$, using $n$ i.i.d.~samples $X_1,...,X_n \stackrel{i.i.d.}{\sim} P$ and $Y_1,...,Y_n \stackrel{i.i.d.}{\sim} Q$ from each distribution. Specifically, we assume that $P$ and $Q$ lie in a smaller subspace $\H_b \subseteq \H_a$ parameterized by a $\Z$-indexed real family $b = \{b_z\}_{z \in \Z}$. In this setting, we study the minimax $\L^2$ error $M(a,b)$ of estimating $\langle P, Q \rangle_a$, over $P$ and $Q$ lying in a (unit) ellipsoid with respect to $\|\cdot\|_b$; that is, the quantity
\begin{equation}
M(a,b)
  := \inf_{\hat S} \quad \sup_{\|P\|_b,\|Q\|_b \leq 1} \quad \E_{\substack{X_1,...,X_n \sim P \\ Y_1,...,Y_n \sim Q}} \left[ \left| \hat S(X_1,...,X_n,Y_1,...,Y_n) - \langle P, Q \rangle_a \right|^2 \right],
  \label{eq:minimax}
\end{equation}
where the infimum is taken over all estimators $\hat S$ (i.e., all complex-valued functions $\hat S : \X^{2n} \to \C$ of the data).

We study how the rate of the minimax error $M(a,b)$ is primarily governed by the rates at which $a_z$ and $b_z$ decay to $0$ as $\|z\| \to \infty$.\footnote{By equivalence of finite-dimensional norms, the choice of norm here affects only constant factors.} This has been studied extensively in the Sobolev (or polynomial-decay) case, where, for some $t > s \geq 0$, $a_z = \|z\|^{-s}$ and $b_z = \|z\|^{-t}$, corresponding to estimation of $s$-order Sobolev semi-inner products under $t$-order Sobolev smoothness assumptions on the Lebesgue density functions $p$ and $q$ of $P$ and $Q$
% \Bk{do you mean $p$ and $q$ to be the densities associated with $P$ and $Q$} \CC{Yes, clarified.}
(as described in Example~\ref{ex:Sobolev} below) \citep{bickel88squaredDerivatives,donoho90minimaxQuadratic,laurent00adaptiveQuadratic,singh16sobolev}. In this case, the rate of $M(a,b)$ has been identified (by \citet{bickel88squaredDerivatives}, \citet{donoho90minimaxQuadratic}, and \citet{singh16sobolev}, in increasing generality) as
\footnote{Here and elsewhere, $\asymp$ denotes equality up to constant factors.}
\begin{equation}
M(a,b) \asymp \max \left\{ n\inv, n^{-\frac{8(t - s)}{4t + D}} \right\},
\label{eq:Sobolev_minimax_rate}
\end{equation}
so that the ``parametric'' rate $n\inv$ dominates when $t \geq 2s + D/4$, and the slower rate $n^{-\frac{8(t - s)}{4t + D}}$ dominates otherwise. \citet{laurent00adaptiveQuadratic} additionally showed that, for $t < 2s + D/4$, $M(a,b)$ increases by a factor of $\left( \log n \right)^{\frac{4(t - s)}{4t + D}}$ in the ``adaptive'' case, when the tail index $t$ is not assumed to be known to the estimator.

However, the behavior of $M(a,b)$ for other (non-polynomial) decay rates of $a$ and $b$ has not been studied, despite the fact that, as discussed in Section~\ref{sec:examples}, other rates of decay of $a$ and $b$, such as Gaussian or exponential decay, correspond to inner products and assumptions commonly considered in nonparametric statistics.
The goal of this paper is therefore to understand the behavior of $M(a,b)$ for general sequences $a$ and $b$.

Although our results apply more generally, to simply summarize our results, consider the case where $a$ and $b$ are ``radial''; i.e. $a_z$ and $b_z$ are both functions of some norm $\|z\|$.
Under mild assumptions, we show that the minimax convergence rate is then a function of the quantities
\[A_{\zeta_n} = \sum_{\|z\| \leq \zeta_n} a_z^{-2}
  \quad \text{ and } \quad B_{\zeta_n} = \sum_{\|z\| \leq \zeta_n} b_z^{-2},\]
which can be thought of as measures of the ``strengths'' of $\|\cdot\|_a$ and $\|\cdot\|_b$, for a particular choice of a ``smoothing'' (or ``truncation'') parameter $\zeta_n \in (0, \infty)$. Specifically, we show
\begin{equation}
  M(a,b) \asymp \max \left\{ \left( \frac{A_{\zeta_n}}{B_{\zeta_n}} \right)^2, \frac{1}{n} \right\},
  \quad \text{ where } \quad
  \zeta_n^D n^2 = B_{\zeta_n}^2.
\label{eq:simplified_minimax_rate}
\end{equation}
While \eqref{eq:simplified_minimax_rate} is difficult to simplify or express in a closed form in general, it is quite simple to compute given the forms of $a$ and $b$. In this sense, \eqref{eq:simplified_minimax_rate} might be considered as an analogue of the Le Cam equation~\citep{yang1999minimax} (which gives a similar implicit formula for the minimax rate of nonparametric density estimation in terms of covering numbers) for estimating inner products and related quantities. It is easy to check that, in the Sobolev case (where $a_z = \|z\|^{-s}$ and $b_z = \|z\|^{-t}$ decay polynomially), \eqref{eq:simplified_minimax_rate} recovers the previously known rate~\eqref{eq:Sobolev_minimax_rate}. Moreover, our assumptions are also satisfied by other rates of interest, such as exponential (where $a_z = e^{-s\|z\|_1}$ and $b_z = e^{-t\|z\|_1}$) and Gaussian (where $a_z = e^{-s\|z\|_2^2}$ and $b_z = e^{-t\|z\|_2^2}$) rates, for which we are the first to identify minimax rates. As in the Sobolev case, the rates here exhibit the so-called ``elbow'' phenomenon, where the convergence rates is ``parametric'' (i.e., of order $\asymp 1/n$) when $t$ is sufficiently large relative to $s$, and slower otherwise. However, for rapidly decaying $b$ such as in the exponential case, the location of this elbow no longer depends directly on the dimension $D$; the parametric rate is achieved as soon as $t \geq 2s$.

We note that, in all of the above cases, the minimax rate \eqref{eq:simplified_minimax_rate} is achieved by a simple bilinear estimator:
\[\hat S_{\zeta_n}
  := \sum_{\|z\| \leq \zeta_n} \frac{\hat \phi_P(z) \overline{\hat \phi_Q(z)}}{a_z^2},\]
where
\[\hat \phi_P(z) := \frac{1}{n} \sum_{i = 1}^n \psi_z(X_i)
  \quad \text{ and } \quad
  \hat \phi_Q(z) := \frac{1}{n} \sum_{i = 1}^n \psi_z(Y_i)\]
are linear estimates of $\phi_P(z)$ and $\phi_Q(z)$, and $\zeta_n \geq 0$ is a tuning parameter.
% In the radial case, we can take $Z$ to be a ball $Z = \{z \in \Z : \|z\| \leq \zeta_n\}$ for some $\zeta_n \in (0, \infty)$, and $\hat S_Z$ is then computable in $O(n \zeta_n^D \log \zeta_n)$ time via the fast Fourier transform; since the optimal choice of $\zeta_n$ often scales as $O(n^{1/D})$ or $O(\log n)$, this can be reasonably efficient in practice.
We also show that, in many cases, a rate-optimal $\zeta_n$ can be chosen adaptively (i.e., without knowledge of the space $\H_b$ in which $P$ and $Q$ lie).
% The most difficult question is then to choose the tuning parameter $Z$, which relies knowledge of the sequence $b$. To construct an adaptive estimator, a simple penalization approach, along the lines of that proposed by \citet{laurent00adaptiveQuadratic}, may suffice. However, given that minimax rates for adaptive estimators are known to differ in some cases from those of non-adaptive estimators under known smoothness, it is likely that a significantly different analysis will be required for this type of estimator, and we leave this for future work.

\subsection{Motivating Examples}
\label{sec:examples}
Here, we briefly present some examples of products $\langle \cdot, \cdot \rangle_a$ and spaces $\H_a$ of the form \eqref{eq:estimand} that are commonly encountered in statistical theory and functional analysis. In the following examples, the base measure on $\X$ is taken to be the Lebesgue measure $\mu$, and ``probability densities'' are with respect to $\mu$. Also, for any integrable function $f \in \L^1(\X)$, we use $\tilde f_z = \int_\X f \psi_z \, d\mu$ to denote the $z^{th}$ Fourier coefficient of $f$ (where $\psi_z$ is the $z^{th}$ Fourier basis element as in~\eqref{def:fourier_coefficients}).

The simplest example is the standard $\L^2$ inner product:
\begin{example}
In the ``unweighted'' case where $a_z = 1$ for all $z \in \Z$, $\H_a$ includes the usual space $\L^2(\X)$ of square-integrable probability densities on $\X$, and, for $P$ and $Q$ with square-integrable densities $p,q \in \L^2(\X)$, we have
\[\langle p, q \rangle_a = \int_\X p(x) q(x) \, dx.\]
\end{example}

Typically, however, we are interested in weight sequences such that $a_z \to 0$ as $\|z\| \to \infty$ and $\H_a$ will be strictly smaller than $\L^2(\X)$ to ensure that $\langle \cdot, \cdot \rangle_a$ is finite-valued; this corresponds intuitively to requiring additional smoothness of functions in $\H$. Here are two examples widely used in statistics:

\begin{example}
If $\H_K$ is a reproducing kernel Hilbert space (RKHS) with a symmetric, translation-invariant kernel $K(x, y) = \kappa(x - y)$ (where $\kappa \in \L^2(\X)$), one can show via Bochner's theorem (see, e.g., Theorem 6.6 of \citep{wendland2005scattered}) that the semi-inner product induced by the kernel can be written in the form
\[\langle f, g \rangle_{\H_K}
  := \sum_{z \in \Z} \tilde \kappa_z^{-2} \tilde f_z \overline{\tilde g_z}.\]
Hence, setting each $a_z = \langle \kappa, \psi_z \rangle = \tilde \kappa_z$, $\H_a$ contains any distributions $P$ and $Q$ on $\X$ with densities $p,q \in \H_K = \{p \in \L^2 : \langle p, p \rangle_{\H_K} < \infty\}$, and we have $\langle P, Q \rangle_a = \langle p, q \rangle_{\H_K}$.
\end{example}

\begin{example}
For $s \in \mathbb{N}$, $\H^s$ is the $s$-order Sobolev space
\begin{equation*}
\H^s := \left\{ f \in \L^2(\X) : f \text{ is $s$-times weakly differentiable with } 
%\footnotemark~with } 
f^{(s)} \in \L^2(\X) \right\},
\end{equation*}
%\footnotetext{Weak differentiability is a generalization of differentiability based on integration by parts. See \citet{evans98PDEs} or \citet{leoni09SobolevSpaces} for definitions, examples, and discussion of weak differentiability in Sobolev spaces.}
endowed with the semi-inner product of the form
\begin{equation}
\langle p, q \rangle_{\H^s}
  := \left\langle p^{(s)}, q^{(s)} \right\rangle_{\L^2(\X)}
  = \sum_{z \in \Z} |z|^{2s} \tilde f_z \overline{\tilde g_z}
  \label{eq:Sobolev_inner_product}
\end{equation}
where the last equality follows from Parseval's identity. 
Indeed, \eqref{eq:Sobolev_inner_product} is commonly used to generalize $\langle f, g \rangle_{\H^s}$, for example, to non-integer values of $s$. Thus, setting $a_z = |z|^{-s}$, $\H_a$ contains any distributions $P,Q \in \P$ with densities $p,q \in \H^s$, and, moreover,
we have $\langle P, Q \rangle_a = \langle p, q \rangle_{\H^s}$.
Note that, when $s \geq D/2$, one can show via Bochner's theorem that $\H^s$ is in fact also an RKHS, with symmetric, translation-invariant kernel defined as above by $\kappa(x) = \sum_{z \in \Z} z^{-s} \psi_z$.
\label{ex:Sobolev}
\end{example}

\subsection*{Paper Organization}
The remainder of this paper is organized as follows:
In Section~\ref{sec:notation}, we provide notation needed to formally state our estimation problem, given in Section~\ref{sec:problem_statement}.
Section~\ref{sec:related_work} reviews related work on estimation of functionals of probability densities, as well as some applications of this work.
Sections~\ref{sec:upper_bounds} and \ref{sec:lower_bounds} present our main theoretical results, with upper bounds in Sections~\ref{sec:upper_bounds} and minimax lower bounds in Section~\ref{sec:lower_bounds}; proofs of all results are given in Appendix~\ref{appendix:proofs}.
Section \ref{sec:special_cases} expands upon these general results in a number of important special cases.
Finally, we conclude in Section~\ref{sec:discussion} with a discussion of broader consequences and avenues for future work.

\section{Notation}
\label{sec:notation}
We assume the sample space $\X \subseteq \R^D$ is a compact subset of $\R^D$, and we use $\mu$ to denote the usual Lebesgue measure on $\X$. We use $\{\psi_z\}_{z \in \mathbb{Z}^D}$ to denote the standard orthonormal Fourier basis of $\L_2(\X)$, indexed by $D$-tuples of integer frequencies $z \in \mathbb{Z}^D$.
For any function $f \in \L_2(\X)$ and $z \in \mathbb{Z}^D$, we use
\[\tilde f_z := \int_\X f(x) \overline{\psi_z(x)} \, d\mu(x)\]
to denote the $z^{th}$ Fourier coefficient of $f$ (i.e., the projection of $f$ onto $\psi_z$), and for any probability distribution $P \in \P$, we use the same notation
\[\phi_P(z)
  := \E_{X \sim P} \left[ \overline{\psi_z(X)} \right]
  = \int_\X \overline{\psi_z(x)} dP(x)\]
to denote the characteristic function of $P$.

We will occasionally use the notation $\|z\|$ for indices $z \in \Z$. Due to equivalence of finite dimensional norms, the exact choice of norm affects only constant factors; for concreteness, one may take the Euclidean norm.

For certain applications, it is convenient to consider only a subset $\Z \subseteq \mathbb{Z}^D$ of indices of interest (for example, Sobolev seminorms are indexed only over $\Z = \{z \in \mathbb{Z}^D : z_1,...,z_D \neq 0 \}$). The subset $\Z$ may be considered arbitrary but fixed in our work.

Given two $(0,\infty)$-valued sequences\footnote{A more proper mathematical term for $a$ and $b$ would be \emph{net}.}
$a = \{a_z\}_{z \in \Z}$ and $b = \{b_z\}_{z \in \Z}$, we are interested in products of the form
\[\langle f, g \rangle_a
  := \sum_{z \in \Z} \frac{\tilde f_z \overline{\tilde g_z}}{a_z^2},\]
and their induced (semi)norms $\|f\|_a = \sqrt{\langle f, f \rangle_a}$ over spaces of the form\footnote{Specifically, we are interested in probability densities, which lie in the simplex $\mathcal{P} := \{f \in \L_1(\X) : f \geq 0, \int_\X f \, d\mu = 1\}$, so that we should write, e.g., $p, q \in \H \cap \mathcal{P}$. Henceforth, ``density'' refers to any function lying in $\mathcal{P}$.}
\[\H_a = \left\{ f \in \L_2(\X) : \|f\|_a < \infty \right\}\]
(and similarly when replacing $a$ by $b$).
Typically, we will have $a_z,b_z \to 0$ and $\frac{b_z}{a_z} \to 0$ whenever $\|z\| \to \infty$, implying the inclusion $\H_b \subseteq \H_a \subseteq \L^2(\X)$.

\section{Formal Problem Statement}
\label{sec:problem_statement}
Suppose we observe $n$ i.i.d. samples $X_1,...,X_n \stackrel{i.i.d.}{\sim} P$ and $n$ i.i.d. samples $Y_1,...,Y_n \stackrel{i.i.d.}{\sim} Q$, where $P$ and $Q$ are (unknown) distributions lying in the (known) space $\H_a$. We are interested in the problem of estimating the inner product \eqref{eq:estimand}, along with the closely related (squared) seminorm and pseudometric given by
\begin{equation}
\|P\|_a^2
  := \langle P, P \rangle_a
  \quad \text{ and } \quad
\|P - Q\|_a^2
  := \|P\|_a^2 + \|Q\|_a^2 - 2 \langle P, Q \rangle_a.
\label{eq:norms_and_distances}
\end{equation}

We assume $P$ and $Q$ lie in a (known) smaller space $\H_b \subseteq \H_a$, and we are specifically interested in identifying, up to constant factors, the minimax mean squared (i.e., $\L^2$) error $M(a,b)$ of estimating $\langle P, Q \rangle_a$ over $P$ and $Q$ lying in a unit ellipsoid with respect to $\|\cdot\|_b$; that is, the quantity
\begin{equation}
M(a,b)
  := \inf_{\hat S} \quad \sup_{\|P\|_b,\|Q\|_b \leq 1} \quad \E_{\substack{X_1,...,X_n \sim p, \\ Y_1,...,Y_n \sim q}} \left[ \left| \hat S - \langle p, q \rangle_a \right|^2 \right],
  \label{eq:minimax}
\end{equation}
where the infimum is taken over all estimators (i.e., all functions $\hat S : \R^{2n} \to \C$ of the data $X_1,...,X_n,Y_1,...,Y_n$).

\section{Related Work}
\label{sec:related_work}
%\Bk{Usually, it is not an accepted notion to start a subsection right after a section. So, please include a sentence or two mentioning what you are going to do in this section before moving onto subsections.} \CC{I've added a bit below.}
This section reviews previous studies on special cases of the problem we study, as well as work on estimating related functionals of probability distributions, and a few potential applications of this work in statistics and machine learning.

\subsection{Prior work on special cases}
While there has been substantial work on estimating unweighted $\L_2$ norms and distances of densities~\citep{schweder75L2,anderson94L2TwoSampleTest,gine08simpleL2}, to the best of our knowledge, most work on the more general problem of estimating \emph{weighted} inner products or norms has been on estimating Sobolev quantities (see Example~\ref{ex:Sobolev} in Section~\ref{sec:intro}) by \citet{bickel88squaredDerivatives}, \citet{donoho90minimaxQuadratic}, and \citet{singh16sobolev}.
\citet{bickel88squaredDerivatives} considered the case of integer-order Sobolev norms,
which have the form
\begin{equation}
\|f\|_{\H^s}^2 = \|f^{(s)}\|_{\L^2(\X)}^2 = \int \left( f^{(s)}(x) \right)^2 \, dx,
\label{eq:classical_integer_Sobolev}
\end{equation}
for which they upper bounded the error of an estimator based on plugging a kernel density estimate into \eqref{eq:classical_integer_Sobolev} and then applying an analytic bias correction.
They also derived matching minimax lower bounds for this problem.\footnote{\citet{bickel88squaredDerivatives} actually make H\"older assumptions on their densities (essentially, an $\L_\infty$ bound on the derivatives of the density), rather than our slightly milder Sobolev assumption (essentially, an $\L_2$ bound on the derivative). However, as we note in Section~\ref{sec:discussion}, these assumptions are closely related such that the results are comparable up to constant factors.}
\citet{singh16sobolev} proved rate-matching upper bounds on the error of a much simpler inner product estimator (generalizing an estimator proposed by \citet{donoho90minimaxQuadratic}), which applies for arbitrary $s \in \R$.
Our upper and lower bounds are strict generalizations of these results.
Specifically, relative to this previous work on the Sobolev case, our work makes advances in three directions:
\begin{enumerate}
\item
We consider estimating a broader class of inner product functionals $\langle p, q \rangle_z$, for arbitrary sequences $\{a_z\}_{z \in \Z}$. The Sobolev case corresponds to $a_z = \|z\|^{-s}$ for some $s > 0$.
\item
We consider a broader range of assumptions on the true data densities, of the form $\|p\|_b,\|q\|_b < \infty$, for arbitrary sequences $\{b_z\}_{z \in \Z}$. The Sobolev case corresponds to $b_z = \|z\|^{-t}$ for some $t > 0$.
\item
We prove lower bounds that match our upper bounds, thereby identifying minimax rates. For many cases, such as Gaussian or exponential RKHS inner products or densities, these results are the first concerning minimax rates, and, even in the Sobolev case, our lower bounds address some previously open cases (namely, non-integer $s$ and $t$, and $D > 1$.
\end{enumerate}

The closely related work of \citet{fan91quadratic} also generalized the estimator of \citet{donoho90minimaxQuadratic}, and proved (both upper and lower) bounds on $M(a,b)$ for somewhat more general sequences, and also considered norms with exponent $p \neq 2$ (i.e., norms not generated by an inner product, such as those underlying a broad class of Besov spaces). However, his analysis placed several restrictions on the rates of $a$ and $b$; for example, it requires
\[\sup_{Z \subseteq \Z} \frac{|Z| \sup_{z \in Z} a_z^{-2}}{\sum_{z \in Z} a_z^{-2}} < \infty
  \quad \text{ and } \quad
  \sup_{Z \subseteq \Z} \frac{|Z| \sup_{z \in Z} b_z^{-2}}{\sum_{z \in Z} b_z^{-2}} < \infty.\]
This holds when $a$ and $b$ decay polynomially, but fails in many of the cases we consider, such as exponential decay. The estimation of norms with $p \neq 2$ and $a$ and $b$ decaying non-polynomially, therefore, remains an important unstudied case, which we leave for future work.

Finally, we note that, except \citet{singh16sobolev}, all the above works have considered only $D = 1$ (i.e., when the sample space $\X \subseteq \R$), despite the fact that $D$ can play an important role in the convergence rates of the estimators. The results in this paper hold for arbitrary $D \geq 1$.

\subsection{Estimation of related functionals}

There has been quite a large amount of recent work~\citep{nguyen2010divergenceFunctionals,liu2012exponential,moon14divergencesEnsemble,singh14RenyiDivergence,singh14densityFunctionals,krishnamurthy14RenyiAndFriends,moon14divergencesConfidence,krishnamurthy15L2Divergence,kandasamy15vonMises,gao15stronglyDependent,gao15localGaussian,mukherjee15lepski,mukherjee16adaptive,moon16improvingConvergence,singh16kNN,berrett16kNNentropy,gao17density,gao2017estimating,jiao2017nearest,han2017optimal,noshad2017direct,wisler2017direct,singh2017nonparanormal,noshad2018scalable,bulinski2018statistical}
% \textcolor{red}{[TODO: Any more references here?]}
on practical estimation of nonlinear integral functionals of probability densities, of the form
% \footnote{Some of this work has considered more general functionals that can depend on multiple densities or their derivatives.}
\begin{equation}
F(p) = \int_\X \varphi(p(x)) \, dx,
\label{eq:integral_functional}
\end{equation}
where $\phi : [0, \infty) \to \R$ is nonlinear but smooth.
Whereas minimax optimal estimators have been long established, their computational complexity typically scales as poorly as $O(n^3)$~\citep{birge95integralFunctionals,laurent96integralFunctionals,kandasamy15vonMises}. Hence, this recent work has focused on analyzing more computationally efficient (but less statistically efficient) estimators, as well as on estimating information-theoretic quantities such as variants of entropy, mutual information, and divergence, for which $\phi$ can be locally non-smooth (e.g., $\phi = \log$), and can hence follow somewhat different minimax rates.

As discussed in detail by \citet{laurent96integralFunctionals}, under Sobolev smoothness assumptions on $p$, estimation of quadratic functionals (such as those considered in this paper) is key to constructing minimax rate-optimal estimators for general functionals of the form~\eqref{eq:integral_functional}.
The reason for this is that minimax rate-optimal estimators of $F(p)$ can often be constructed by approximating a second-order Taylor (a.k.a., von Mises~\citep{kandasamy15vonMises}) expansion of $F$ around a density estimate $\hat p$ of $p$ that is itself minimax rate-optimal (with respect to integrated mean squared error). Informally, if we expand $F(p)$ as
\begin{equation}
F(p)
  = F(\hat p)
  + \langle \nabla F(\hat p), p - \hat p \rangle_{\L^2}
  + \left\langle p - \hat p, (\nabla^2 F(\hat p)) p - \hat p \right\rangle_{\L^2}
  + O \left( \|p - q\|_{\L^2}^3 \right),
  \label{eq:von_Mises_expansion}
\end{equation}
where $\nabla F(\hat p)$ and $\nabla^2 F(\hat p)$ are the first and second order Frechet derivatives of $F$ at $\hat p$.
% \Bk{I am somewhat confused with the notation in the second order term. Should it be $\langle p - \hat p, \nabla^2 F(\hat p) (p - \hat p)\rangle_{\L^2}$}
In the expansion~\eqref{eq:von_Mises_expansion}, the first term is a simple plug-in estimate, and the second term is linear in $p$, and can therefore be estimated easily by an empirical mean. The remaining term is precisely a quadratic functional of the density, of the type we seek to estimate in this paper.
Indeed, to the best of our knowledge, this is the approach taken by \emph{all} estimators that are known to achieve minimax rates \citep{birge95integralFunctionals,laurent96integralFunctionals,krishnamurthy14RenyiAndFriends,kandasamy15vonMises,mukherjee15lepski,mukherjee16adaptive} for general functionals of the form~\eqref{eq:integral_functional}.

Interestingly, the estimators studied in the recent papers above are all based on either
kernel density estimators~\citep{singh14RenyiDivergence,singh14densityFunctionals,krishnamurthy14RenyiAndFriends,krishnamurthy15L2Divergence,kandasamy15vonMises,moon16improvingConvergence,mukherjee15lepski,mukherjee16adaptive}
or $k$-nearest neighbor methods~\citep{moon14divergencesEnsemble,moon14divergencesConfidence,singh16kNN,berrett16kNNentropy,gao17density}. This contrasts with our approach, which is more comparable to orthogonal series density estimation; given the relative efficiency of computing orthogonal series estimates (e.g., via the fast Fourier transform), it may be desirable to try to adapt our estimators to these classes of functionals.

When moving beyond Sobolev assumptions, only estimation of very specific functionals has been studied. For example, under RKHS assumptions, only estimation of maximum mean discrepancy (MMD)\citep{gretton12kernel,ramdas15decreasingPower,tolstikhin16minimaxMMD}, has received much attention.
Hence, our work significantly expands our understanding of minimax functional estimation in this setting. More generally, our work begins to provide a framework for a unified understanding of functional estimation across different types of smoothness assumptions.

Along a different line, there has also been some work on estimating $\L^p$ norms for \emph{regression functions}, under similar Sobolev smoothness assumptions~\citep{lepski99regressionNorm}.
However, the problem of norm estimation for regression functions turns out to have quite different statistical properties and requires significantly different estimators and analysis, compared to norm estimation for density functions.
Generally, the problem for densities is statistically easier in terms of having a faster convergence rate under a comparable smoothness assumption; this is most obvious when $p = 1$, since the $\L^1$ norm of a density is always $1$, while the $\L^1$ norm of a regression function is less trivial to estimate.
However, this is true more generally as well. For example, \citet{lepski99regressionNorm} showed that, under $s$-order Sobolev assumptions, the minimax rate for estimating the $\L^2$ norm of a $1$-dimensional regression function (up to $\log$ factors) is $\asymp n^{-\frac{4s}{4s + 1}}$, whereas the corresponding rate for estimating the $\L^2$ norm of a density function is $\asymp n^{-\min\left\{\frac{8s}{4s + 1}, 1 \right\}}$, which is parametric when $s \geq 1/4$. To the best of our knowledge, there has been no work on the natural question of estimating Sobolev or other more general quadratic functionals of regression functions.

\subsection{Applications}

Finally, although this paper focuses on estimation of general inner products from the perspective of statistical theory, we mention a few of the many applications that motivate the study of this problem.

Estimates of quadratic functionals can be directly used for nonparametric goodness-of-fit, independence, and two-sample testing~\citep{anderson94L2TwoSampleTest,dumbgen98goodnessOfFit,ingster12nonparametric,goria05new,pardo05statistical,chwialkowski15fastTwoSample}. They can also by used to construct confidence sets for a variety of nonparametric objects~\citep{li89honestConfidence,baraud04confidence,genovese05waveletRegression}, as well as for parameter estimation in semi-parametric models~\citep{wolsztynski05minimumEntropyEstimation}.

In machine learning, Sobolev-weighted distances can also be used in transfer learning~\citep{du17hypothesisTransferLearning}
and transduction learning~\citep{quadrianto09transduction}
to measure relatedness between source and target domains, helping to identify when transfer can benefit learning.
Semi-inner products can be used as kernels over probability distributions, enabling generalization of a wide variety of statistical learning methods from finite-dimensional vectorial inputs to nonparametric distributional inputs~\citep{sutherland16thesis}.
This \emph{distributional learning} approach has been applied to many diverse problems, including image classification~\citep{Poczos:2011:NDE:3020548.3020618,poczos12imageClassification}, galaxy mass estimation~\citep{ntampaka15galaxyMass}, ecological inference~\citep{flaxman15election,flaxman16election}, aerosol prediction in climate science~\citep{szabo15learningTheoryDistributions}, and causal inference~\citep{lopez15causalTheory}.
Finally, it has recently been shown that the losses minimized in certain implicit generative models can be approximated by Sobolev and related distances~\citep{liang2017well}.
Further applications of these quantities can be found in~\citep{principe10information}.

\section{Upper Bounds}
\label{sec:upper_bounds}
In this section, we provide upper bounds on minimax risk.
Specifically, we propose estimators for semi-inner products, semi-norms, and pseudo-metrics, and bound the risk of the semi-inner product estimator; identical bounds (up to constant factors) follow easily for semi-norms and pseudo-metrics.

\subsection{Proposed Estimators}
\label{subsec:proposed_estimator}

Our proposed estimator $\hat S_Z$ of $\langle P, Q \rangle_a$ consists of simply plugging estimates of $\phi_P$ and $\tilde Q$ into a truncated version of the summation in Equation~\eqref{eq:estimand}. Specifically, since
\[\phi_P(z) = \E_{X \sim P} \left[ \overline{\psi_z(X)} \right],\]
% \Bk{In notation section, you had complex conjugation over $\psi_z$ which was absent in introduction and also here. So please keep it uniform.} \CC{Fixed.}
we estimate each $\phi_P(z)$ by $\hat \phi_P(z) := \frac{1}{n} \sum_{i = 1}^n \psi_z(X_i)$ and each $\phi_Q(z)$ by $\hat \phi_Q(z) := \frac{1}{n} \sum_{i = 1}^n \psi_z(Y_i)$. Then, for some finite set $Z \subseteq \Z$ (a tuning parameter to be chosen later) our estimator $\hat S_Z$ for the product~\eqref{eq:estimand} is
\begin{equation}
\hat S_Z := \sum_{z \in Z} \frac{\hat \phi_P(z) \overline{\hat \phi_Q(z)}}{a_z^2}.
\label{eq:inner_product_estimator}
\end{equation}
To estimate the squared semi-norm $\|P\|_a^2$ from a single sample $X_1,\dots,X_n \stackrel{i.i.d.}{\sim} P$, we use
\begin{equation}
\hat N_Z := \sum_{z \in Z} \frac{\hat \phi_P(z) \overline{\hat \phi_P(z)'}}{a_z^2}.
\label{eq:norm_estimator}
\end{equation}
where $\hat \phi_P(z)$ is estimated using the first half $X_1,\dots,X_{\lfloor n/2 \rfloor}$ of the sample,
$\hat \phi_P(z)'$ is estimated using the second half $X_{\lfloor n/2 \rfloor + 1},\dots,X_n$ of the sample.
While it is not clear that sample splitting is optimal in practice, it allows us to directly apply convergence results for the semi-inner product, which assume the samples from the two densities are independent.

To estimate the squared pseudo-metric $\|P - Q\|_a^2$ from two samples $X_1,\dots,X_n \stackrel{i.i.d.}{\sim} P$ and $Y_1,\dots,Y_n \stackrel{i.i.d.}{\sim} Q$, we combine the above inner product and norm estimators according to the formula~\eqref{eq:norms_and_distances}, giving
\[\hat \rho_Z = \hat N_Z + \hat M_Z - 2 \hat S_Z,\]
where $\hat M_Z$ denotes the analogue of the norm estimator~\eqref{eq:norm_estimator} applied to $Y_1,\dots,Y_n$.

\subsection{Bounding the risk of \texorpdfstring{$\hat S_Z$}{SZ}}

Here, we state upper bounds on the bias, variance, and mean squared error of the semi-inner product estimator $\hat S_Z$, beginning with an easy bound on the bias of $\hat S_Z$ (proven in Appendix~\ref{sec:bias_bound}):
\begin{proposition}[Upper bound on bias of $\hat S_Z$]
Suppose $P, Q \in \H_b$. Then,
\begin{equation}
\left| \B \left[ \hat S_Z \right] \right|
%   := \left| \E \left[ \hat S_Z \right] - \langle P, Q \rangle_a \right|
  \leq \|P\|_b \|Q\|_b \sup_{z \in \Z \sminus Z} \frac{b_z^2}{a_z^2},
\label{ineq:bias_bound}
\end{equation}
where $\B \left[ \hat S_Z \right] := \E \left[ \hat S_Z \right] - \langle P, Q \rangle_a$ denotes the bias of $\hat S_Z$.
\label{prop:bias_bound}
\end{proposition}
Note that for the above bound to be non-trivial, we require $b_z \to 0$ faster than $a_z$ as $\|z\| \to \infty$ which ensures that $\sup_{z \in \Z \sminus Z} \frac{b_z}{a_z} < \infty$. While \eqref{ineq:bias_bound} does not explicitly depend on the sample size $n$, in practice, the parameter set $Z$ will be chosen to grow with $n$, and hence the supremum over $\Z \sminus Z$ will decrease monotonically with $n$.
Next, we provide a bound on the variance of $\hat S_Z$, whose proof, given in Appendix~\ref{sec:variance_bound}, is more involved.
\begin{proposition}[Upper bound on variance of $\hat S_Z$]
Suppose $P, Q \in \H_b$. Then,
\begin{equation}
\Var[\hat S_Z]
  \leq \frac{2\|P\|_2 \|Q\|_2}{n^2} \sum_{z \in Z} \frac{1}{a_z^4}
  + \frac{\|Q\|_b^2 \|P\|_b + \|P\|_b^2 \|Q\|_b}{n} R_{a,b,Z}
  + \frac{2\|P\|_a^2 \|Q\|^2_a}{n}
\label{ineq:var_bound}
\end{equation}
where $\mathbb{V}$ denotes the variance operator and
\begin{equation}
R_{a,b,Z} := \left( \sum_{z \in Z} \frac{b_z^4}{a_z^8} \right)^{1/4} \left( \sum_{z \in Z} \left( \frac{b_z}{a_z^2} \right)^8 \right)^{1/8}
\left( \sum_{z \in Z} b_z^8 \right)^{1/8}.
\label{eq:RabZ}
\end{equation}
\label{prop:var_bound}
\end{proposition}

Having bounded the bias and variance of the estimator $\hat S_Z$, we now turn to the mean squared error (MSE). Via the usual decomposition of MSE into (squared) bias and variance, Propositions~\ref{prop:bias_bound} and~\ref{prop:var_bound} together immediately imply the following bound:
\begin{theorem}[Upper bound on MSE of $\hat S_Z$]
Suppose $P, Q \in \H_b$. Then,
\begin{align}
\MSE \left[ \hat S_Z \right]
\notag
& \leq \|P\|_b^2 \|Q\|_b^2 \sup_{z \in \Z \sminus Z} \frac{b_z^4}{a_z^4}
  + \frac{2\|P\|_2 \|Q\|_2}{n^2} \sum_{z \in Z} \frac{1}{a_z^4} \\
\label{ineq:MSE_bound}
& \quad + \frac{\|Q\|_b^2 \|P\|_b + \|P\|_b^2 \|Q\|_b}{n} R_{a,b,Z}
  + \frac{2\|P\|^2_a \|Q\|^2_a}{n},
\end{align}
where $R_{a,b,Z}$ is as defined in \eqref{eq:RabZ}.
\label{thm:MSE_bound}
\end{theorem}
\begin{corollary}[Norm estimation]
In the particular case of norm estimation (i.e., when $Q = P$), this simplifies to:
\begin{equation}
\MSE \left[ \hat S_Z \right]
  \leq \|P\|_b^4 \sup_{z \in \Z \sminus Z} \frac{b_z^4}{a_z^4}
  + \frac{2\|P\|_2^2}{n^2} \sum_{z \in Z} \frac{1}{a_z^4}
  + \frac{2\|P\|_b^3}{n} R_{a,b,Z}
  + \frac{2\|P\|_a^4}{n}.
\label{ineq:norm_MSE_bound}
\end{equation}
\label{corr:norm_MSE_bound}
\end{corollary}

\subsection{Discussion of Upper Bounds}
Two things might stand out that distinguish the above variance bound from many other nonparametric variance bounds:
First, the rate depends on the smoothness of $P, Q \in \H_b$. Smoothness assumptions in nonparametric statistics are usually needed only to bound the bias of estimators~\citep{Tsybakov:2008:INE:1522486}.
The reason the smoothness appears in this variance bound is that the estimand in Equation~\eqref{eq:estimand} includes products of the Fourier coefficients of $P$ and $Q$. Hence, the estimates $\hat \phi_P(z)$ of $\phi_P(z)$ are scaled by $\hat \phi_Q(z)$, and vice versa, and as a result, the decay rates of $\phi_P(z)$ and $\phi_Q(z)$ affect the variance of the tails of $\hat S_Z$. One consequence of this is that the convergence rates exhibit a phase transition, with a parametric convergence rate when the tails of $\phi_P$ and $\tilde Q$ are sufficiently light, and a slower rate otherwise.

Second, the bounds are specific to the Fourier basis (as opposed to, say, any uniformly bounded basis, e.g., one with $\sup_{z \in \Z, x \in \X} |\psi_z(x)| \leq 1$). The reason for this is that, when expanded, the variance includes terms of the form $\E_{X \sim P}[\phi_y(X)\psi_z(X)]$, for some $y \neq z \in \Z$. In general, these covariance-like terms are difficult to bound tightly; for example, the uniform boundedness assumption above would only give a bound of the form $\E_{X \sim P}[|\phi_y(X)\psi_z(X)|] \leq \min\{\phi_P(y), \phi_P(z)\}$. For the Fourier basis, however, the recurrence relation $\phi_y \psi_z = \phi_{y + z}$ allows us to bound $\E_{X \sim P}[\phi_y(X)\psi_z(X)] = \E_{X \sim P}[\phi_{y + z}(X)] = \phi_P(y + z)$ in terms of assumptions on the decay rates of the coefficients of $P$. It turns out that $\phi_P(y + z)$ decays significantly faster than $\min\{\phi_P(y), \phi_P(z)\}$, and this tighter bound is needed to prove optimal convergence rates.

More broadly, this suggests that convergence rates for estimating inner products in terms of weights in a particular basis may depend on algebraic properties of that basis. For example, another common basis, the Haar wavelet basis, satisfies a different recurrence relation: $\phi_y \psi_z \in \{0, \phi_y, \psi_z\}$, depending on whether (and how) the supports of $\phi_y$ and $\psi_z$ are nested or disjoint. We leave investigation of this and other bases for future work.

Clearly, $\sup_{Z \subseteq \Z} R_{a,b,Z} < \infty$ if and only if $b_z^4/a_z^8$ is summable (i.e., $\sum_{z \in \Z} b_z^4/a_z^8 < \infty$). Thus, assuming $|\Z| = \infty$, this already identifies the precise condition required for the minimax rate to be parametric. When
% For the optimal choice $Z^*$ of the parameter $Z$
% \[Z^* := \argmin_{Z \subseteq \Z} \sup_{z \in \Z \sminus Z} \frac{b_z^4}{a_z^4}
%   + \frac{1}{n^2} \sum_{z \in Z} \frac{1}{a_z^4},\]
it is the case that
\[\frac{R_{a,b,Z}}{n} \in O \left( \sup_{z \in \Z \sminus Z} \frac{b_z^4}{a_z^4} + \frac{1}{n^2} \sum_{z \in Z} \frac{1}{a_z^4} \right),\]
the third term in~\eqref{ineq:MSE_bound} will be dominated by the first and third terms, and so
% \Bk{I would agree that if the above assumption of $R_{a,b,Z^*}$ being lower order holds, then the following claims are all ok. Somehow, the way the above sentence is written implies that the lower order of $R_{a,b,Z^*}$ implicitly holds. Is that true? If so, I dont see why.}
% \CC{
% \[R_{a,b,Z} := \left( \sum_{z \in Z} \frac{b_z^4}{a_z^8} \right)^{1/4} \left( \sum_{z \in Z} \left( \frac{b_z}{a_z^2} \right)^8 \right)^{1/8}
% \left( \sum_{z \in Z} b_z^8 \right)^{1/8}.\]
% If $b_z^4/a_z^8$ is summable, then $R_{a,b,Z}/n \in O(n\inv)$, and so this term is negligible. Otherwise, $b_z^4/a_z^8 \in \Omega(z^{-D})$.
% Suffices to show:
% \[\frac{1}{n} \left( \sum_{z \in Z} \frac{b_z^4}{a_z^8} \right)^{1/2}
%   \leq \sup_{z \in \Z \sminus Z} \frac{b_z^2}{a_z^2} \frac{1}{n} \left( \sum_{z \in Z} a_z^{-4} \right)^{1/2}\]
% \[\left( \sum_{z \in Z} \frac{b_z^4}{a_z^8} \right)^{1/2}
%   \leq \sup_{z \in \Z \sminus Z} \frac{b_z^2}{a_z^2} \left( \sum_{z \in Z} a_z^{-4} \right)^{1/2}\]
% \begin{align*}
% R_{a,b,Z} := \left( \sum_{z \in Z} \frac{b_z^4}{a_z^8} \right)^{1/4} \left( \sum_{z \in Z} \frac{b_z^8}{a_z^{16}} \right)^{1/8}
% \left( \sum_{z \in Z} b_z^8 \right)^{1/8}
% \end{align*}
% }
% Therefore, 
the upper bound simplifies to order
\begin{equation}
\MSE \left[ \hat S_{Z^*} \right] \lesssim \frac{1}{n} + \min_{Z \subseteq \Z} \left[\sup_{z \in \Z \sminus Z} \frac{b_z^4}{a_z^4}
  + \frac{1}{n^2} \sum_{z \in Z} \frac{1}{a_z^4}\right].
  \label{eq:simplified_upper_bound}
\end{equation}
This happens for every choice of $a_z$ and $b_z$ we consider in this paper, including the Sobolev (polynomial decay) case and the RKHS case. However, simplifying the bound further requires some knowledge of the form of $a$ and/or $b$, and we develop this in several cases in Section~\ref{sec:special_cases}. In Section~\ref{sec:discussion}, we also consider some heuristics for approximately simplifying \eqref{eq:simplified_upper_bound} in certain settings.

\section{Lower Bounds}
\label{sec:lower_bounds}
In this section, we provide a lower bound on the minimax risk of the estimation problems described in Section~\ref{sec:problem_statement}.
Specifically, we use a standard information theoretic framework to lower bound the minimax risk for semi-norm estimation; bounds of the same rate follow easily for inner products and pseudo-metrics. In a wide range of cases, our lower bound matches the MSE upper bound (Theorem~\ref{thm:MSE_bound}) presented in the previous section.

\begin{theorem}[Lower Bound on Minimax MSE]
Suppose $\X$ has finite base measure $\mu(\X) = 1$ and suppose the basis $\{\psi_z\}_{z \in \Z}$ contains the constant function $\phi_0 = 1$ and is uniformly bounded (i.e., $\sup_{z \in \Z} \|\psi_z\|_\infty < \infty$). Define $Z_{\zeta_n} := \{z \in \mathbb{Z}^D : \|z\|_\infty \leq \zeta_n\}$, $A_{\zeta_n} := \sum_{z \in Z_{\zeta_n}} a_z^{-2}$ and $B_{\zeta_n} := \sum_{z \in Z_{\zeta_n}} b_z^{-2}$. If $B_{\zeta_n} \in \Omega \left( \zeta_n^{2D} \right)$, then we have the minimax lower bound
\[\inf_{\hat S} \sup_{\|h\|_b \leq 1} \E \left[ \left( \hat S - \|h\|_a^2 \right)^2 \right]
  \in \Omega \left( \max \left\{ \frac{A_{\zeta_n}^2}{B_{\zeta_n}^2}, n\inv \right\} \right),\]
where $\zeta_n$ is chosen to satisfy $B_{\zeta_n}^2 \asymp \zeta_n^D n^2$.
Also, if $B_{\zeta_n} \in o \left( \zeta_n^{2D} \right)$, then we have the (looser) minimax lower bound
\[\inf_{\hat S} \sup_{\|h\|_b \leq 1} \E \left[ \left( \hat S - \|h\|_a^2 \right)^2 \right]
  \in \Omega \left( \max \left\{ \frac{A_{n^{2/3D}}^2}{n^{4/3}}, n\inv \right\} \right).\]
\label{thm:lower_bound}
\end{theorem}

\begin{remark}
The uniform boundedness assumption permits the Fourier basis, our main case of interest, but also allows other bases (see, e.g., the ``generalized Fourier bases'' used in Corollary 2.2 of \citet{liang2017well}).
\end{remark}

\begin{remark}
The condition that $B_{\zeta_n} \in \Omega \left( \zeta_n^{2D} \right)$ is needed to ensure that the ``worst-case'' densities we construct in the proof of Theorem~\ref{thm:lower_bound} are indeed valid probability densities (specifically, that they are non-negative). Hence, this condition would no longer be necessary if we proved results in the simpler Gaussian sequence model, as in many previous works on this problem (e.g., \citep{cai99adaptive,cai05nonquadraticQuadratic}). However, when $B_{\zeta_n} \in o \left( \zeta_n^{2D} \right)$, density estimation, and hence the related problem of norm estimation, become asymptotically easier than the analogous problems under the Gaussian sequence model.
\end{remark}

\begin{remark}
Intuitively, the ratio $A_{\zeta_n}/B_{\zeta_n}$ measures the relative strengths of the norms $\|\cdot\|_a$ and $\|\cdot\|_b$. As expected, consistent estimation is possible if and only if $\|\cdot\|_b$ is a stronger norm than $\|\cdot\|_a$.
\end{remark}

\section{Special Cases}
\label{sec:special_cases}
In this section, we develop our lower and upper results for several special cases of interest. The results of this section are summarized in Table~\ref{tab:minimax_rates}.\vspace{1mm}\\
{\bf Notation:} Here, for simplicity, we assume that the estimator $\hat S_Z$ uses a choice of $Z$ that is symmetric across dimensions; in particular,
$Z = \prod_{j = 1}^D \{\phi_{-\zeta_n},...,\phi_0,...,\phi_{\zeta_n}\}$ (for some $\zeta_n \in \N$ depending on $n$) is the Cartesian product of $D$ sets of the first $2\zeta_n + 1$ integers.
Throughout this section, we use $\lesssim$ and $\gtrsim$ to denote inequality up to $\log n$ factors. Although we do not explicitly discuss estimation of $\L^2$ norms, it appears as a special case of the Sobolev case with $s = 0$.
% \Bk{In the following, you have used $\zeta_n$ at some places and $\zeta$ at some other. Please keep it consistent.} \CC{I've replaced $\zeta$ with $\zeta_n$ everywhere.}
%\begin{enumerate}
%\item
%{\bf Sobolev Case:} 
\subsection{Sobolev}
For some $s, t \geq 0$, $a_z = \|z\|^{-s}$ and $b_z = \|z\|^{-t}$.\vspace{1mm}\\
{\bf Upper Bound:} By Proposition~\ref{prop:bias_bound}, $\B[\hat S_Z] \lesssim \zeta_n^{2(s - t)}$, and, by Proposition~\ref{prop:var_bound},
\[\Var[\hat S_Z]
%   \lesssim \frac{\zeta_n^{4s + D}}{n^2}
%   + \frac{\zeta_n^{2s - t + D/4} \zeta_n^{2s - t} \zeta_n^{-t + D/4}}{n}
%   + \frac{1}{n}
  \lesssim \frac{\zeta_n^{4s + D}}{n^2}
  + \frac{\zeta_n^{4s - 3t + D/2}}{n}
  + \frac{1}{n}\]
Thus,
\[\MSE[\hat S_Z]
  \lesssim \zeta_n^{4(s - t)} + \frac{\zeta_n^{4s + D}}{n^2} + \frac{\zeta_n^{4s - 3t + D/2}}{n} + \frac{1}{n}.\]
One can check that $\zeta_n^{4(s - t)} + \frac{\zeta_n^{4s + D}}{n^2}$ is minimized when $\zeta_n \asymp n^{\frac{2}{4t + D}}$, and that, for this choice of $\zeta_n$, the $\frac{\zeta_n^{4s - 3t + D/2}}{n}$ term is of lower order, giving the convergence rate
\[\MSE[\hat S_Z] \asymp n^{\frac{8(s - t)}{4t + D}}.\]
{\bf Lower Bound:} Note that $A_{\zeta_n} \asymp \zeta_n^{2s + D}$ and $B_{\zeta_n} \asymp \zeta_n^{2t + D}$. Solving $B_{\zeta_n}^2 = \zeta_n^D n^2$ gives $\zeta_n \asymp n^{\frac{2}{4t + D}}$. Thus, Theorem~\ref{thm:lower_bound} gives a minimax lower bound of
\[\inf_{\hat S} \sup_{\|p\|_b, \|q\|_b \leq 1} \E \left[\left( \hat S - \langle p, q \rangle_b \right)^2 \right]
  \gtrsim \frac{A_{\zeta_n}^2}{B_{\zeta_n}^2}
  = \zeta_n^{4(s - t)}
  = n^{\frac{8(s - t)}{4t + D}},\]
matching the upper bound. Note that the rate is parametric ($\asymp n\inv$) when $t \geq 2s + D/4$, and slower otherwise.
%\item
%{\bf Gaussian RKHS:}
\subsection{Gaussian RKHS}
For some $t \geq s \geq 0$, $a_z = e^{-s\|z\|_2^2}$ and $b_z = e^{-t\|z\|_2^2}$.\vspace{1mm}\\
{\bf Upper Bound:} By Proposition~\ref{prop:bias_bound}, $\B[\hat S_Z] \lesssim e^{2(s - t)\zeta_n^2}$.
If we use the upper bound
\[\sum_{z \in Z} e^{\theta \|z\|_2^2}
  \leq C_\theta \zeta_n^D e^{\theta \zeta_n^2},\]
for any $\theta > 0$ and some $C_\theta > 0$, then Proposition~\ref{prop:var_bound} gives
\begin{align*}
\Var[\hat S_Z]
& = \frac{\zeta_n^D e^{4s\zeta_n^2}}{n^2}
  + \frac{\zeta_n^{89D/20} e^{(4s - 3t)\zeta_n^2}}{n}
  + \frac{1}{n}.
\end{align*}
Thus,
\[\MSE[\hat S_Z]
  \lesssim e^{4(s - t)\zeta_n^2}
  + \frac{\zeta_n^D e^{s\zeta_n^2}}{n^2}
  + \frac{\zeta_n^{89D/20} e^{(4s - 3t)\zeta_n^2}}{n}
  + \frac{1}{n}.\]
One can check that $e^{4(s - t)\zeta_n^2} + \frac{\zeta_n^D e^{s\zeta_n^2}}{n^2}$ is minimized when $\zeta_n \asymp \sqrt{\frac{\log n}{2t}}$, and that, for this choice of $\zeta_n$, the $\zeta_n^{89D/20} e^{(4s - 3t)\zeta_n^2}$ term is of lower order, giving an MSE convergence rate of
\[\MSE[\hat S_Z]
  \lesssim n^{\frac{2(s - t)}{t}}
  = n^{2(s/t - 1)}.\]
{\bf Lower Bound:} Again, we use the bound
\[A_{\zeta_n}
  = \sum_{z \in Z_{\zeta_n}} e^{2s\|z\|_2^2}
  \lesssim \zeta_n^D e^{2s\zeta_n^2},\]
as well as the trivial lower bound $B_{\zeta_n}
  = \sum_{z \in Z_{\zeta_n}} e^{2s\|z\|_2^2}
  \geq e^{2s\zeta_n^2}$.
Solving $B_{\zeta_n}^2 = \zeta_n^D n^2$ gives $\zeta_n \asymp \sqrt{\frac{\log n}{2t}}$ up to $\log \log n$ factors. Thus, ignoring $\log n$ factors, Theorem~\ref{thm:lower_bound} gives a minimax lower bound of
\[\inf_{\hat S} \sup_{\|p\|_b, \|q\|_b \leq 1} \E \left[\left( \hat S - \langle p, q \rangle_b \right)^2 \right] \gtrsim n^{\frac{2(s - t)}{t}},\]
for some $C > 0$, matching the upper bound rate. Note that the rate is parametric when $t \geq 2s$, and slower otherwise.
%\item
%{\bf Exponential RKHS:}
\subsection{Exponential RKHS}
For some $t \geq s \geq 0$, $a_z = e^{-s\|z\|_1}$ and $b_z = e^{-t\|z\|_1}$.\vspace{1mm}\\
{\bf Upper Bound:} By Proposition~\ref{prop:bias_bound}, $\B[\hat S_Z] \lesssim e^{2(s - t) \zeta_n}$. Since, for fixed $D$,
\[\sum_{z \in Z} e^{r\|z\|_1}
  \asymp e^{r\zeta_n + D}
  \asymp e^{r\zeta_n},\]
by Proposition~\ref{prop:var_bound}, we have
\begin{align*}
\Var[\hat S_Z]
& \asymp \frac{e^{4s\zeta_n}}{n^2}
  + \frac{e^{(4s - 3t)\zeta_n}}{n}
  + \frac{1}{n},
\end{align*}
giving a mean squared error bound of
\[\MSE[\hat S_Z]
  \asymp e^{4(s - t)\zeta_n} + \frac{e^{4s\zeta_n}}{n^2} + \frac{e^{(4s - 3t) \zeta_n}}{n} + \frac{1}{n}.\]One can check that $e^{4(s - t)\zeta_n} + \frac{e^{4s\zeta_n}}{n^2}$ is minimized when $\zeta_n \asymp \frac{\log n}{2t}$, and that, for this choice of $\zeta_n$, the $\frac{e^{(4s - 3t) \zeta_n}}{n}$ term is of lower order, giving an MSE convergence rate of
\[\MSE[\hat S_Z]
  \lesssim n^{\frac{2(s - t)}{t}}
  = n^{2(s/t - 1)}.\]
{\bf Lower Bound:} Note that $A_{\zeta_n} \asymp e^{2s\zeta_n}$ and $B_{\zeta_n} = e^{2t\zeta_n}$. Solving $B_{\zeta_n}^2 = \zeta_n^D n^2$ gives, up to $\log\log n$ factors, $\zeta_n \asymp \frac{\log n}{2t}$. Thus, Theorem~\ref{thm:lower_bound} gives a minimax lower bound of
\[\inf_{\hat S} \sup_{\|p\|_b, \|q\|_b \leq 1} \E \left[\left( \hat S - \langle p, q \rangle_b \right)^2 \right] \gtrsim n^{\frac{2(s - t)}{t}},\]
for some $C > 0$, matching the upper bound rate. Note that the rate is parametric when $t \geq 2s$, and slower otherwise.
%\item
%{\bf Logarithmic decay:} 
\subsection{Logarithmic decay}
For some $t \geq s \geq 0$, $a_z = \left( \log \|z\| \right)^{-s}$ and $b_z = \left( \log \|z\| \right)^{-t}$. Note that, since our lower bound requires $B_{\zeta_n} \in \Omega(\zeta_n^{2D})$, we will only study the upper bound for this case.\vspace{1mm}\\
{\bf Upper Bound:} By Proposition~\ref{prop:bias_bound}, $\B[\hat S_Z] \lesssim \left( \log \zeta_n \right)^{2(s - t)}$. By the upper bound
\[\sum_{z \in Z_{\zeta_n}} \left( \log \|z\| \right)^\theta
  \leq C_\theta \zeta_n^D \left( \log \zeta_n \right)^\theta,\]
for any $\theta > 0$ and some $C_\theta > 0$, Proposition~\ref{prop:var_bound} gives
\begin{align*}
\Var[\hat S_Z]
& = \frac{\zeta_n^D \left( \log \zeta_n \right)^{4s}}{n^2}
  + \frac{\zeta_n^{89D/20} \left( \log \zeta_n \right)^{4s - 3t}}{n}
  + \frac{1}{n},
\end{align*}
giving a mean squared error bound of
\[\MSE[\hat S_Z]
  \lesssim \left( \log \zeta_n \right)^{4(s - t)}
  + \frac{\zeta_n^D \left( \log \zeta_n \right)^{4s}}{n^2}
  + \frac{\zeta_n^{89D/20} \left( \log \zeta_n \right)^{4s - 3t}}{n}
  + \frac{1}{n}.\]
One can check that $\left( \log \zeta_n \right)^{4(s - t)} + \frac{\zeta_n^{89D/20} \left( \log \zeta_n \right)^{4s - 3t}}{n}$ is minimized when $\zeta_n^{89D/20} \left( \log \zeta_n \right)^{4t + D} \asymp n$, and one can check that, for this choice of $\zeta_n$, the $\frac{\zeta_n^D \left( \log \zeta_n \right)^{4s}}{n^2}$ term is of lower order.

Thus, up to $\log n$ factors, $\zeta_n \asymp n^{2/D}$, and so, up to $\log \log n$ factors,
\[\left( \log \zeta_n \right)^{4(s - t)}
  \asymp \left( \log n \right)^{4(s - t)}.\]
%\item
%{\bf Sinc RKHS:}
\subsection{Sinc RKHS} 
For any $s \in (0,\infty)^D$, the $\sinc_s$ kernel, defined by
\[K_{\sinc}^s(x, y) =
\prod_{j = 1}^d \frac{s_j}{\pi} \sinc \left( \frac{x_j - y_j}{s_j} \right),\]
where
\[\sinc(x) =
  \left\{
    \begin{array}{ll}
      \frac{\sin(x)}{x} & \text{ if } x \neq y \\
      1 & \text{ else }
  \end{array}
  \right.
,\]
generates the RKHS $\H_{\sinc}^s = \left\{ f \in \L^2 : \|f\|_{K_{\sinc}^s} < \infty \right\}$,
of band-limited functions, where the norm is generated by the inner product
$\langle f, g \rangle_{K_{\sinc}^s} = \langle f, g \rangle_a$,
where $a_z = 1_{\{|z| \leq s\}}$ (with the convention that $\frac{0}{0} = 0$). If we assume that $p \in \H_{\sinc}^t$, where $t \leq s$, then fixing $Z := \{z \in \mathbb{Z}^D : |z| \leq s \}$, by Proposition~\ref{prop:bias_bound}, $\B[\hat S_Z] = 0$, and, by Proposition~\ref{prop:var_bound}, one can easily check that $\Var[\hat S_Z] \lesssim n\inv$. Thus, without any assumptions on $P$, we can always estimate $\|P\|_{K_{\sinc}^s}$ at the parametric rate.
% \Bk{I do not follow the symbol $\preceq$. I assumed it is a element-wise comparison but then it does not make sense in the definition of the set $Z$.} \CC{I'm not sure why I used $\preceq$ here. I've changed it to $\leq$.}
%\end{enumerate}

\begin{table}
\begin{center}
\begin{tabular}{l|c|c|c|c}
                                       & $b_z = \log^{-t} \|z\|$ & $b_z = \|z\|^{-t}$            & $b_z = e^{-t\|z\|_1}$ & $b_z = e^{-t\|z\|_2^2}$ \\
\hline
$a_z = \log^{-s} \|z\|$ & $\lesssim\left(\log n\right)^{4(s-t)}$ & $\max \left\{ n\inv, n^{\frac{-8t}{4t + D}} \right\}$                       & $n\inv$               & $n\inv$              \\
\hline
$a_z = \|z\|^{-s}$                     & $\infty$                               & $\max \left\{ n\inv, n^{\frac{8(s - t)}{4t + D}} \right\}$ & $n\inv$               & $n\inv$                 \\
\hline
$a_z = e^{-s\|z\|_1}$                  & $\infty$                               & $\infty$                      & $\max \left\{ n\inv, n^{2\frac{s - t}{t}} \right\}$      & $n\inv$                 \\
\hline
$a_z = e^{-s\|z\|_2^2}$                & $\infty$                               & $\infty$                      & $\infty$              & $\max \left\{ n\inv, n^{2\frac{s - t}{t}} \right\}$
\end{tabular}
\end{center}
\caption{Minimax convergence rates for different combinations of $a_z$ and $b_z$. Results are given up to $\log n$ factors, except the case when both $a_z$ and $b_z$ are logarithmic, which is given up to $\log \log n$ factors. Note that, in this last case, only the upper bound is known. A value of $\infty$ indicates that the estimand itself may be $\infty$ and consistent estimation is impossible.
}
\label{tab:minimax_rates}
\end{table}

\section{Discussion}
\label{sec:discussion}

In this paper, we focused on the case of inner product weights and density coefficients in the Fourier basis, which play well-understood roles in widely used spaces such as Sobolev spaces and reproducing kernel Hilbert spaces with translation-invariant kernels.

For nearly all choices of weights $\{a_z\}_{z \in \Z}$ and $\{b_z\}_{z \in \Z}$, ignoring the parametric $1/n$ term that appears in both the upper and lower bounds, the upper bound boils down to
\[\min_{\zeta_n \in \N} \frac{b_{\zeta_n}^4}{a_{\zeta_n}^4} + \frac{\sum_{z \in Z_{\zeta_n}} a_z^{-4}}{n^2},\]
or, equivalently,
\[\frac{b_{\zeta_n}^4}{a_{\zeta_n}^4}
  \quad \text{ where } \quad
  \frac{b_{\zeta_n}^4}{a_{\zeta_n}^4}
  = \frac{\sum_{z \in Z_{\zeta_n}} a_z^{-4}}{n^2}\]
and the lower bound boils down to
\[\left( \frac{\sum_{z \in Z_{\zeta_n}} a_z^{-2}}{\sum_{z \in Z_{\zeta_n}} b_z^{-2}} \right)^2,
  \quad \text{ where } \quad
  \left( \sum_{z \in Z_{\zeta_n}} b_z^{-2} \right)^2 = \zeta_n^D n^2.\]
These rates match if
\begin{equation}
\frac{a_{\zeta_n}^{-4}}{b_{\zeta_n}^{-4}}
  \asymp \left( \frac{\sum_{z \in Z_{\zeta_n}} a_z^{-2}}{\sum_{z \in Z_{\zeta_n}} b_z^{-2}} \right)^2
\quad \text{ and } \quad \frac{b_{\zeta_n}^4 \left( \sum_{z \in Z_{\zeta_n}} b_z^{-2} \right)^2}{a_{\zeta_n}^4 \sum_{z \in Z_{\zeta_n}} a_z^{-4}}
  \asymp \zeta_n^D
\label{eq:rates_match}
\end{equation}
Furthermore, if the equations in \eqref{eq:rates_match} hold modulo logarithmic factors, then the upper and lower bounds match modulo logarithmic factors. This holds almost automatically if $b_z$ decays exponentially or faster, since, then, $\zeta_n$ grows logarithmically with $n$. Noting that the lower bound requires $B_{\zeta_n} \in \Omega \left( \zeta_n^{2D} \right)$, this also holds automatically if $b_z = |z|^t$ with $t \geq D/2$.

Table~\ref{tab:minimax_rates} collects the derived minimax rates for various standard choices of $a$ and $b$. For entries below the diagonal, $b_z/a_z \to \infty$ as $\|z\| \to \infty$, and so $\H_b \not\subseteq \H_a$. As a result, consistent estimation is not possible in the worst case. The diagonal entries of Table~\ref{tab:minimax_rates}, for which $a$ and $b$ have the same form, are derived in Section~\ref{sec:special_cases} directly from our upper and lower bounds on $M(a,b)$. These cases exhibit a phase transition, with convergence rates depending on the parameters $s$ and $t$. When $t$ is sufficiently larger than $s$, the variance is dominated by the low-order terms of the estimand~\eqref{eq:estimand}, giving a convergence rate of $\asymp n\inv$. Otherwise, the variance is dominated by the tail terms of~\ref{eq:estimand}, in which case minimax rates depend smoothly on $s$ and $t$. This manifests in the $\max \{n\inv, n^{R(s,t)} \}$ form of the minimax rates, where $R$ is non-decreasing in $s$ and non-increasing in $t$.

Notably, the data dimension $D$ plays a direct role in the minimax rate only in the Sobolev case when $t < 2s + D/4$. Otherwise, the role of $D$ is captured entirely within the assumption that $p,q \in \H_b$. This is consistent with known rates for estimating other functionals of densities under strong smoothness assumptions such as the RKHS assumption~\citep{gretton12kernel,ramdas15decreasingPower}.

Finally, we note some consequences for more general (non-Hilbert) Sobolev spaces $W^{s,p}$, defined for $s \geq 0$, $p \geq 1$ as the set of functions in $\L^p$ having weak $s^{th}$ derivatives in $\L^p$. The most prominent example is that of the H\"older spaces $W^{s,\infty}$ of essentially bounded functions having essentially bounded $s^{th}$ weak derivatives; H\"older spaces are used widely in nonparametric statistics~\citep{bickel88squaredDerivatives,Tsybakov:2008:INE:1522486}.
Recall that, for $p \leq q$ and any $s \geq 0$, these spaces satisfy the embedding $\W^{s,q} \subseteq \W^{s,p}$~\citep{villani1985another}, and that $W^{s,2} = \H^s$. Then, for $P, Q \in \W^{t,p}$ our upper bound in Theorem~\ref{thm:MSE_bound} implies an identical upper bound when $p \geq 2$, and our lower bound in Theorem~\ref{thm:lower_bound} implies an identical lower bound when $p \leq 2$.

Further work is needed to verify tightness of these bounds for $p \neq 2$.
Moreover, while this paper focused on the Fourier basis, it is also interesting to consider other bases, which may be informative in other spaces. For example, wavelet bases are more natural representations in a wide range of Besov spaces~\citep{donoho95wavelet}. It is also of interest to consider non-quadratic functionals as well as non-quadratic function classes. In these cases simple quadratic estimators such as those considered here may not achieve the minimax rate, but it may be possible to correct this with simple procedures such as thresholding, as done, for example, by \citet{cai05nonquadraticQuadratic} in the case of $\L_p$ balls with $p < 2$. Finally, the estimators considered here require some knowledge of the function class in which the true density lies. It is currently unclear whether and how the various strategies for designing adaptive estimators, such as block-thresholding~\citep{cai99adaptive} or Lepski's method~\citep{lepski97optimal}, which have been applied to estimate quadratic functionals over $\L_p$ balls and Besov spaces~\citep{efromovich96optimalAdaptive,cai06adaptive}, may confer adaptivity when estimating functionals over general quadratically weighted spaces.

%\Bk{I wouldnt write a separate conclusion section. I will combine this with the discussion section. You can start the discussion with the way you started in conclusion and slowly move onto discussing your results as you did and then explain the exploratory work. If you do this (which I would suggest), please make changes to the paper organization section also.} \CC{I've made this change.}

% \begin{supplement}
% \sname{Supplement A}\label{suppA}
% \stitle{Proofs Title of the Supplement A}
% \slink[url]{http://www.e-publications.org/ims/support/dowload/imsart-ims.zip}
% \sdescription{Dum esset rex in
% accubitu suo, nardus mea dedit odorem suavitatis. Quoniam confortavit
% seras portarum tuarum, benedixit filiis tuis in te. Qui posuit fines tuos}
% \end{supplement}
\section{Proofs}
\label{appendix:proofs}
In this section, we present the proofs of main results.
\subsection{Proof of Proposition~\ref{prop:bias_bound}}
\label{sec:bias_bound}
We first bound the bias $\left| \E \left[ \hat S_Z \right] - \langle P, Q \rangle_a \right|$, where randomness is over the data $X_1,...,X_n \stackrel{i.i.d.}{\sim} P$, and $Y_1,...,Y_n \stackrel{i.i.d.}{\sim} Q$. Since
\[\hat S_Z = \sum_{z \in Z} \frac{\hat \phi_P(z) \overline{\hat \phi_Q(z)}}{a_z^2}\] is bilinear in $\hat P$ and $\hat Q$, which are independent, and
\[\hat \phi_P(z) = \frac{1}{n} \sum_{i = 1}^n \psi_z(X_i)
  \quad \text{ and } \quad
  \hat \phi_Q(z) = \frac{1}{n} \sum_{i = 1}^n \psi_z(Y_i)\] are unbiased estimators of $\phi_P(z) = \E_{X \sim P} \left[ \psi_z(X) \right]$ and $\phi_Q(z) = \E_{Y \sim Q} \left[ \psi_z(Y) \right]$, respectively, we have that
\[\E \left[ \hat S_Z \right]
 = \sum_{z \in Z} \frac{\E\left[\hat \phi_P(z) \right] \overline{\E\left[\hat \phi_Q(z) \right]}}{a_z^2}
 = \sum_{z \in Z} \frac{\phi_P(z) \overline{\phi_Q(z)}}{a_z^2}.\]
Hence, the bias is
\[\left| \E \left[ \hat S_Z \right] - \langle P, Q \rangle_a \right|
 = \left| \sum_{z \in \Z \sminus Z} \frac{\phi_P(z) \overline{\phi_Q(z)}}{a_z^2} \right|.\]
If $p, q \in \H_b \subset \H_a \subseteq \L^2$ defined by 
\[\H_b
 := \left\{
      f \in \L^2 : \|f\|_b := \sum_{z \in \Z} \frac{\tilde f_z^2}{b_z^2} < \infty
    \right\},\]
then, applying Cauchy-Schwarz followed by H\"older's inequality, we have
% twice (or, equivalently, the generalized H\"older's inequality with powers $2$, $2$, and $\infty$),
\begin{align*}
\left| \E \left[ \hat S_Z \right] - \langle P, Q \rangle_a \right|
 = \left| \sum_{z \in \Z \sminus Z} \frac{\phi_P(z) \overline{\phi_Q(z)}}{a_z^2} \right|
& \leq \sqrt{\sum_{z \in \Z \sminus Z} \frac{|\phi_P(z)|^2}{a_z^2} \sum_{z \in \Z \sminus Z} \frac{|\phi_Q(z)|^2}{a_z^2}} \\
& = \sqrt{\sum_{z \in \Z \sminus Z} \frac{b_z^2}{a_z^2} \frac{|\phi_P(z)|^2}{b_z^2} \sum_{z \in \Z \sminus Z} \frac{b_z^2}{a_z^2} \frac{|\phi_Q(z)|^2}{b_z^2}} \\
& \leq \|P\|_b \|Q\|_b \sup_{z \in \Z \sminus Z} \frac{b_z^2}{a_z^2}.
\end{align*}

Note that this recovers the bias bound of \citet{singh16sobolev} in the Sobolev case: If 
$a_z = z^{-s}$ and $b_z = z^{-t}$ with $t \geq s$, then
\[\left| \E \left[ \hat S_Z \right] - \langle P, Q \rangle_a \right|
  \leq \|P\|_b \|Q\|_b |Z|^{2(s - t)},\]
where $|Z|$ denotes the cardinality of the index set $Z$.

\subsection{Proof of Proposition~\ref{prop:var_bound}}
\label{sec:variance_bound}
In this section, we bound the variance of $\Var[S_Z]$, where, again, randomness is over the data $X_1,...,X_n,Y_1,...,Y_n$. The setup and first several steps of our proof are quite general, applying to arbitrary bases. However, without additional assumptions, our approach eventually hits a roadblock. Thus, to help motivate our assumptions and proof approach, we begin by explaining this general setup in Section~\ref{subsubsec:general_setup}, and then proceed with steps specific to the Fourier basis in Section~\ref{subsubsec:fourier_basis}.

\subsubsection{General Proof Setup}
\label{subsubsec:general_setup}
Our bound is based on the Efron-Stein inequality \citep{efron81jackknife}. For this, suppose that we draw extra independent samples $X_1' \sim p$ and $Y_1' \sim q$, and let $\hat S_Z'$ and $\hat S_Z''$ denote the estimator given in Equation~\eqref{eq:inner_product_estimator} when we replace $X_1$ with $X_1'$ and when we replace $Y_1$ with $Y_1'$, respectively, that is
\begin{equation*}
\hat S_Z' := \sum_{z \in Z} \frac{\hat \phi_P(z)'\overline{\hat \phi_Q(z)}}{a_j^2}
\quad \text{ and } \quad
\hat S_Z'' := \sum_{z \in Z} \frac{\hat \phi_P(z) \overline{\hat \phi_Q(z)'}}{a_j^2},
% \label{eq:resampled_inner_product_estimator}
\end{equation*}
where
\[\hat \phi_P(z)' := \frac{1}{n} \left( \psi_z(X_1') + \sum_{i = 2}^n \psi_z(X_i) \right),
  \text{ and } \quad
  \hat \phi_Q(z)' := \frac{1}{n} \left( \psi_z(Y_1') + \sum_{i = 2}^n \psi_z(Y_i') \right).\]

Then, since $X_1,...,X_n$ and $Y_1,...,Y_n$ are each i.i.d., the Efron-Stein inequality \citep{efron81jackknife} gives
\begin{equation}
\Var\left[ \hat S_Z \right]
  \leq \frac{n}{2} \left( \E \left[ \left| \hat S_Z - \hat S_Z' \right|^2 \right] + \E \left[ \left| \hat S_Z - \hat S_Z'' \right|^2 \right] \right).
\label{ineq:efron_stein}
\end{equation}
We now study just the first term as the analysis of the second is essentially identical. Expanding the definitions of $\hat S_Z$ and $\hat S_Z'$, and leveraging the fact that all terms in $\hat \phi_P(z) - \hat \phi_P(z)'$ not containing $X_1$ or $X_1'$ cancel,
\begin{align}
\E \left[ \left| \hat S_Z - \hat S_Z' \right|^2 \right]
& = \E \left[ \left| \sum_{z \in Z} \frac{\left( \hat \phi_P(z) - \hat \phi_P(z)' \right)\overline{\hat \phi_Q(z)}}{a_z^2} \right|^2 \right] \nonumber\\
& = \E \left[ \sum_{z \in Z} \sum_{w \in Z} \left( \frac{\left( \hat \phi_P(z) - \hat \phi_P(z)' \right)\overline{\hat \phi_Q(z)}}{a_z^2} \right) \left( \frac{\overline{\left( \hat P_w - \hat P_w' \right)}\hat Q_w}{a_w^2} \right) \right] \nonumber\\
& = \frac{1}{n^2} \sum_{z \in Z} \sum_{w \in Z} \E \left[ \hat Q_w \overline{\hat \phi_Q(z)} \frac{\left( \psi_z(X_1) - \psi_z(X_1') \right)\overline{\left( \phi_w(X_1) - \phi_w(X_1') \right)}}{a^2_z a^2_w} \right] \nonumber\\
& = \frac{1}{n^2} \sum_{z \in Z} \sum_{w \in Z} \E \left[ \overline{\hat \phi_Q(z)} \hat Q_w \right] \frac{\E \left[ \left( \psi_z(X_1) - \psi_z(X_1') \right)\overline{\left( \phi_w(X_1) - \phi_w(X_1') \right)} \right]}{a^2_z a^2_w} \nonumber\\
& = \frac{2}{n^2} \sum_{z \in Z} \sum_{w \in Z} \E \left[ \overline{\hat \phi_Q(z)} \hat Q_w \right] \frac{\E \left[ \psi_z(X)\overline{\phi_w(X)} \right] - \E \left[ \psi_z(X) \right] \E \left[ \overline{\phi_w(X)} \right]}{a^2_z a^2_w} \nonumber\\
& = \frac{2}{n^2} \sum_{z \in Z} \sum_{w \in Z} \E \left[ \overline{\hat \phi_Q(z)} \hat Q_w \right] \frac{\E \left[ \psi_z(X)\overline{\phi_w(X)} \right] - \phi_P(z) \overline{\tilde P_w}}{a^2_z a^2_w}.\label{Eq:1}
\end{align}

Expanding the $\E \left[ \overline{\hat \phi_Q(z)} \hat Q_w \right]$ term, we have
\begin{align*}
\E \left[ \overline{\hat \phi_Q(z)} \hat Q_w \right]
& = \frac{1}{n^2} \sum_{i = 1}^n \sum_{j = 1}^n \E \left[ \overline{\psi_z(Y_i)} \phi_w(Y_j) \right] \\
& = \frac{1}{n^2} \sum_{i = 1}^n \E \left[ \overline{\psi_z(Y_i)} \phi_w(Y_i) \right]
  + \frac{1}{n^2} \sum_{i = 1}^n \sum_{j \neq i} \E \left[ \overline{\psi_z(Y_i)} \right] \E \left[ \phi_w(Y_j) \right] \\
& = \frac{1}{n} \E \left[ \overline{\psi_z(Y)} \phi_w(Y) \right]
  + \frac{n - 1}{n} \overline{\phi_Q(z)} \tilde Q_w,
\end{align*}
which combined with Equation \eqref{Eq:1} yields
\begin{eqnarray}
%\notag
\lefteqn{\E \left[ \left| \hat S_Z - \hat S_Z' \right|^2 \right]
 % = \frac{2}{n^2} \sum_{z \in Z} \sum_{w \in Z} \E \left[ \overline{\hat \phi_Q(z)} \hat Q_w \right] \frac{\E \left[ \psi_z(X)\overline{\phi_w(X)} \right] - \phi_P(z) \overline{\tilde P_w}}{a^2_z a^2_w} \\
 = \frac{2}{n^3} \sum_{z \in Z} \sum_{w \in Z} \left( \E \left[ \overline{\psi_z(Y)} \phi_w(Y) \right]
  + (n - 1) \overline{\phi_Q(z)} \tilde Q_w \right)}\nonumber\\
&\qquad\qquad\qquad\qquad\qquad\qquad\qquad \displaystyle\times\, \frac{\E \left[ \psi_z(X)\overline{\phi_w(X)} \right] - \phi_P(z) \overline{\tilde P_w}}{a^2_z a^2_w}. 
\label{eq:general_expression}
\end{eqnarray}
To proceed beyond this point, it is necessary to better understand the covariance-like term $\E \left[ \psi_z(X) \overline{\phi_w(X)} \right]$ appearing in the above equation. If $\{\psi_z\}_{z \in \Z}$ is an arbitrary orthonormal basis, it is difficult to argue more than that, via Cauchy-Schwarz,
\[\left| \E \left[ \psi_z(X)\overline{\phi_w(X)} \right] \right|
  \leq \sqrt{\left| \E \left[ |\psi_z(X)|^2 \right] \E \left[ |\phi_w(X)|^2 \right] \right|}.\]
However, considering, for example, the very well-behaved case when $P$ is the uniform density on $\X$, we would have (since $\{\psi_z\}_{z \in \Z}$ is orthonormal) $\E \left[ |\psi_z(X)|^2 \right] = \E \left[ |\phi_w(X)|^2 \right] = \frac{1}{\mu(\X)}$, which does not decay as $\|z\|, \|w\| \to \infty$. If we were to follow this approach, the Efron-Stein inequality would eventually give a variance bound on $\hat S_Z$ that includes a term of the form
\[\frac{(n - 1)}{n^2 \mu(\X)} \sum_{z, w \in Z} \frac{\overline{\phi_Q(z)} \tilde Q_w}{a_z^2 a_w^2}
  = \frac{(n - 1)}{n^2 \mu(\X)} \left( \sum_{z \in Z} \frac{\phi_Q(z)}{a_z^2} \right)^2
  \leq \frac{1}{n \mu(\X)} \|q\|_b^2 \sum_{z \in Z} \frac{b_z^2}{a_z^4}.\]
While relatively general, this bound is loose, at least in the Fourier case.
Hence, we proceed along tighter analysis that is specific to the Fourier basis.

\subsubsection{Variance Bounds in the Fourier Basis}
\label{subsubsec:fourier_basis}
In the case that $\{\psi_z\}_{z \in \Z}$ is the Fourier basis, the identities $\overline{\psi_z} = \phi_{-z}$ and $\psi_z\phi_w = \psi_{z + w}$ imply that $\E[\psi_z(X)\overline{\phi_w(X)}] = \tilde P_{z - w}$ and $\E[\phi_w(Y)\overline{\psi_z(Y)}] = \tilde Q_{w - z}$, thus, the expression~\eqref{eq:general_expression} simplifies to

\begin{align*}
& \E \left[ \left| \hat S_Z - \hat S_Z' \right|^2 \right]
= \frac{2}{n^3} \sum_{z \in Z} \sum_{w \in Z} \left( \tilde Q_{w - z} + (n - 1) \tilde Q_{-z} \tilde Q_w \right) \frac{\tilde P_{z - w} - \phi_P(z) \tilde {P}_{-w}}{a^2_z a^2_w} \\
& = \frac{2}{n^3} \sum_{z \in Z} \sum_{w \in Z} \frac{\tilde Q_{w - z} \tilde P_{z - w} - \tilde Q_{w - z} \phi_P(z) \tilde P_{-w} + (n - 1) \tilde Q_{-z} \tilde Q_w \tilde P_{z - w} - (n - 1) \tilde Q_{-z} \tilde Q_w \phi_P(z) \tilde P_{-w}}{a^2_z a^2_w}.
\end{align*}
This contains four terms to bound, but they are dominated by the following three main terms:
\begin{equation}
\frac{2}{n^3}
\sum_{z \in Z} \sum_{w \in Z} \frac{\tilde Q_{w - z} \tilde P_{z - w}}{a_z^2 a_w^2},
\label{exp:term1}
\end{equation}
\begin{equation}
\frac{2(n - 1)}{n^3}
\sum_{z \in Z} \sum_{w \in Z} \frac{\overline{\phi_Q(z)} \tilde Q_w \tilde P_{z - w}}{a_z^2 a_w^2},
\label{exp:term2}
\end{equation}
and
\begin{equation}
\frac{2(n - 1)}{n^3}
\sum_{z \in Z} \sum_{w \in Z} \frac{\overline{\phi_Q(z)} \tilde Q_w \phi_P(z) \overline{\tilde P_w}}{a_z^2 a_w^2}.
\label{exp:term3}
\end{equation}
%\vspace{2mm}\\
{\bf Bounding \eqref{exp:term1}:} Applying the change of variables $k = z - w$ gives
\begin{eqnarray}
\left|\frac{2}{n^3} \sum_{z \in Z} \sum_{w \in Z} \frac{\tilde P_{z - w} \tilde Q_{w - z}}{a_z^2a_w^2}\right|&{}={}&\frac{2}{n^3} \left|\sum_{k \in z - Z}\tilde P_k \tilde Q_{-k} \sum_{z \in Z} \frac{1}{a_z^2a^2_{z-k}}\right|
  \stackrel{(*)}{\leq} \frac{2}{n^3} \sum_{k \in z - Z} \left|\tilde P_k \tilde Q_{-k}\right| \sum_{z \in Z} \frac{1}{a_z^4}\nonumber\\
 &{} \leq{}& \frac{2}{n^3} \sum_{k \in \Z} \left|\tilde P_k \tilde Q_{-k}\right| \sum_{z \in Z} \frac{1}{a_z^4}
  \leq \frac{2\|P\|_2\|Q\|_2}{n^3} \sum_{z \in Z} \frac{1}{a_z^4},
  \end{eqnarray}
  where in $(*)$, we use the fact that 
\[f_Z(k) := \sum_{z \in Z} \frac{1}{a_z^2a_{z - k}^2}\]
is the convolution (over $\mathbb{Z}$) of $\{a_z^{-2}\}_{z \in Z}$ with itself, which is always maximized when $k = 0$.\vspace{2mm}\\ 
{\bf Bounding \eqref{exp:term2}:} Applying Cauchy-Schwarz inequality twice, yields
\begin{align}
\left| \sum_{z \in Z} \sum_{w \in Z} \frac{\overline{\phi_Q(z)} \tilde Q_w \tilde P_{z - w}}{a_z^2 a_w^2} \right|
\notag
& = \left| \sum_{z \in Z} \frac{\overline{\phi_Q(z)}}{b_z} \sum_{w \in Z} \frac{b_z}{a_z^2} \frac{\tilde Q_w \tilde P_{z - w}}{a_w^2} \right| \\
\notag
& \leq \|Q\|_b \left( \sum_{z \in Z} \left( \sum_{w \in Z} \frac{b_z}{a_z^2} \frac{\tilde Q_w \tilde P_{z - w}}{a_w^2} \right)^2 \right)^{1/2} \\
\notag
& = \|Q\|_b \left( \sum_{z \in Z} \frac{b_z^2}{a_z^4} \left( \sum_{w \in Z} \frac{\tilde Q_w \tilde P_{z - w}}{a_w^2} \right)^2 \right)^{1/2} \\
\label{ineq:bounding_term2}
& \leq \|Q\|_b \left( \sum_{z \in Z} \frac{b_z^4}{a_z^8} \right)^{1/4} \left( \sum_{z \in Z} \left( \sum_{w \in Z} \frac{\tilde Q_w \tilde P_{z - w}}{a_w^2} \right)^4 \right)^{1/4}.
\end{align}
Note that now we can view the expression
\[\left( \sum_{z \in Z} \left( \sum_{w \in Z} \frac{\tilde Q_w \tilde P_{z - w}}{a_w^2} \right)^4 \right)^{1/4}
  = \left\| \frac{\tilde Q}{a^2} * \tilde P \right\|_4\]
as the $\L_4$ norm of the convolution between the sequence $\tilde Q/a^2$ and the sequence $\tilde P$. To proceed, we apply (a discrete variant of) Young's inequality for convolutions~\citep{beckner75inequalities}, which states that, for constants $\alpha, \beta, \gamma \geq 1$ satisfying $1 + 1/\gamma = 1/\alpha + 1/\beta$ and arbitrary functions $f \in \L^\alpha(\R^D), g \in \L^\beta(\R^D)$,
\[\|f * g\|_\gamma \leq \|f\|_\alpha \|g\|_\beta.\]
Applying Young's inequality for convolutions with powers\footnote{This seemingly arbitrary choice of $\alpha$ and $\beta$ arises from 
analytically minimizing the final bound.} $\alpha = \beta = 8/5$ (so that $\alpha, \beta \geq 1$ and $1/\alpha + 1/\beta = 1 + 1/4$),
gives
\begin{align*}
\left( \sum_{z \in Z} \left( \sum_{w \in Z} \frac{\tilde Q_w \tilde P_{z - w}}{a_w^2} \right)^4 \right)^{1/4}
& \leq \left( \sum_{z \in Z} \frac{\phi_Q(z)^\alpha}{a_z^{2\alpha}} \right)^{1/\alpha} \left( \sum_{z \in Z} \phi_P(z)^\beta \right)^{1/\beta} \\
& = \left( \sum_{z \in Z} \frac{\phi_Q(z)^\alpha}{b_z^\alpha} \frac{b_z^\alpha}{a_z^{2\alpha}} \right)^{1/\alpha}
\left( \sum_{z \in Z} \frac{\phi_P(z)^\beta}{b_z^\beta} b_z^\beta \right)^{1/\beta}.
\end{align*}
Since $2/\alpha = 2/\beta \geq 1$, we can now apply H\"older's inequality to each of the above summations, with powers $\left(2/\alpha, \frac{2\alpha}{2 - \alpha}\right) = \left(2/\beta, \frac{2\beta}{2 - \beta}\right) = \left( \frac{5}{4}, \frac{1}{8} \right)$. This gives
\[\left( \sum_{z \in Z} \frac{\phi_Q(z)^\alpha}{b_z^\alpha} \frac{b_z^\alpha}{a_z^{2\alpha}} \right)^{1/\alpha}
  \leq \|Q\|_b \left( \sum_{z \in Z} \left( \frac{b_z}{a_z^2} \right)^{\frac{2\alpha}{2 - \alpha}} \right)^{\frac{2 - \alpha}{2\alpha}}
  = \|Q\|_b \left( \sum_{z \in Z} \left( \frac{b_z}{a_z^2} \right)^8 \right)^{1/8}\]
and
\[\left( \sum_{z \in Z} \frac{\phi_P(z)^\beta}{b_z^\beta} b_z^\beta \right)^{1/\beta}
  \leq \|P\|_b \left( \sum_{z \in Z} b_z^{\frac{2\beta}{2 - \beta}} \right)^{\frac{2 - \beta}{2\beta}}
  = \|P\|_b \left( \sum_{z \in Z} b_z^8 \right)^{1/8}.\]
Combining these inequalities with inequality~\eqref{ineq:bounding_term2} gives
\[\left| \sum_{z \in Z} \sum_{w \in Z} \frac{\overline{\phi_Q(z)} \tilde Q_w \tilde P_{z - w}}{a_z^2 a_w^2} \right|
  \leq \|Q\|_b^2 \|P\|_b R_{a,b,Z},\]
where $R_{a,b,Z}$ is as in~\eqref{eq:RabZ}.\vspace{2mm}\\
{\bf Bounding \eqref{exp:term3}:} Applying Cauchy-Schwarz yields
\begin{align*}
\frac{2(n - 1)}{n^3} \sum_{z \in Z} \sum_{w \in Z} \frac{\tilde Q_w \overline{\phi_Q(z)} \phi_P(z) \overline{\tilde P_w}}{a_z^2 a_w^2}
& = \frac{2(n - 1)}{n^3} \left( \sum_{z \in Z} \frac{\overline{\phi_Q(z)} \phi_P(z)}{a_z^2} \right) \left( \sum_{w \in Z}\frac{ \tilde Q_w \overline{\tilde P_w}}{a_w^2} \right) \\
& \leq \frac{2}{n^2} \left( \sum_{z \in Z} \frac{|\phi_Q(z)|^2}{a_z^2} \right) \left( \sum_{z \in Z} \frac{|\phi_P(z)|^2}{a_z^2} \right)
  = \frac{2\|P\|_a^2 \|Q\|_a^2}{n^2}.
\end{align*}
Plugging these into Efron-Stein yields the result.
%will eventually give
%\[\Var[\hat S_Z]
%  \leq \frac{2\|P\|_2 \|Q\|_2}{n^2} \sum_{z \in Z} \frac{1}{a_z^4}
%  + \frac{\|Q\|_b^2 \|P\|_b}{n} R_{a,b,Z}
%  + \frac{2\|P\|_a^2 \|Q\|_a^2}{n}.\]

\subsection{Proof of Theorem~\ref{thm:lower_bound}}
\begin{proof}
% \Bk{I would suggest to remove the proof sketch and provide the full proof in a detailed manner. } \CC{I feel the sketch is useful for someone who isn't well-versed with proving minimax lower bounds. But I can remove it, if you think the full proof is clear enough.}\Bk{I think full proof is better. The paper is written for specialists anyway and these techniques are pretty standard. so lets remove the sketch. If you want, you can elaborate the proof more detailedly by highlighting the key ideas.}
% {\bf Sketch:}
The $\Omega(n\inv)$ term of the lower bound, reflecting parametric convergence when the tails of the estimand~\eqref{eq:estimand} are light relative to the first few terms, follows from classic information bounds~\citep{bickel88squaredDerivatives}. We focus on deriving the $\Omega(A_\zeta^2/B_\zeta^2)$ term, reflecting slower convergence when the estimand is dominated by its tail. To do this, we consider the uniform density $\psi_0$ and a family of $2^{|Z_\zeta|}$ small perturbations of the form
\begin{equation}
g_{\zeta,\tau} = \psi_0 + c_\zeta \sum_{z \in Z_\zeta} \tau_z \psi_z,
\label{eq:worst_case_functions}
\end{equation}
where $\zeta \in \N$, $\tau \in \{-1,1\}^{Z_\zeta}$, and $c_\zeta = B_\zeta^{-1/2}$.

We now separately consider the ``smooth'' case, in which $B_\zeta \in \Omega \left( \zeta^{2D} \right)$, and the ``unsmooth'' case, in which $B_\zeta \in o \left( \zeta^{2D} \right)$.

{\bf The smooth case ($B_\zeta \in \Omega \left( \zeta^{2D} \right)$):}
% We now turn to the more involved problem of showing the minimax rate is $\Omega(A_\zeta^2/B_\zeta^2)$, for an appropriate choice of $\zeta$ scaling with $n$, assuming $B_\zeta \in \Omega(\zeta^{2D})$.
By Le Cam's Lemma (see, e.g., Section 2.3 of \citet{Tsybakov:2008:INE:1522486}), it suffices to prove four main claims about the family of $g_{\zeta,\tau}$ functions defined in Equation~\eqref{eq:worst_case_functions}:
\begin{enumerate}
\item
Each $\|g_{\zeta,\tau}\|_b \leq 1$.
\item
Each
\[\inf_{\tau \in \{-1,1\}^{Z_\zeta}} \|g_{\zeta,\tau}\|_a - \|\psi_0\|_a
  \geq \frac{A_\zeta}{B_\zeta}.\]
\item
Each $g_{\zeta,\tau}$ is a density function (i.e., $\int_\X g_{\zeta,\tau} = 1$ and $g_{\zeta,\tau} \geq 0$).
\item
$\zeta$ and $c_\zeta$ are chosen (depending on $n$) such that
\[D_{TV}\left( \psi_0^n, \frac{1}{2^{|Z_\zeta|}} \sum_{\tau \in \{-1,1\}^{Z_\zeta}} g_{\zeta,\tau}^n \right) \leq \frac{1}{2}.\]
\end{enumerate}

For simplicity, for now, suppose $\Z = \N^D$ and $Z_\zeta = [\zeta]^D$. For any $\tau \in \{-1,1\}^{Z_\zeta}$, let
\[g_{\zeta,\tau}
  = \psi_0
  + c_\zeta \sum_{z \in Z_\zeta} \tau_z \psi_z.\]
By setting $c_\zeta = B_\zeta^{-1/2} = \left( \sum_{z \in Z_\zeta} b_z^{-2} \right)^{-1/2}$, we automatically ensure the first two claims:
% [TODO: Clarify how $\psi_0$ term affects norm.]
\[\|g_{\zeta,\tau}\|_b^2
  = c_\zeta^2 \sum_{z \in Z_\zeta} b_z^{-2}
  = 1,\]
and
\[\inf_{\tau \in \{-1,1\}^{Z_\zeta}} \|g_{\zeta,\tau}\|_a^2 - \|\psi_0\|_a^2
  = c_\zeta^2 \sum_{z \in Z_\zeta} a_z^{-2}
  = \frac{A_\zeta}{B_\zeta}.\]

To verify that each $g_{\zeta,\tau}$ is a density, we first note that, since, for $z \neq 0$, $\int_\X \psi_z = 0$, and so
\[\int_\X g_{\zeta,\tau} = \int_\X \psi_0 = 1.\]
Also, since $\psi_0$ is constant and strictly positive and the supremum is taken over all $\tau \in \{-1,1\}^{Z_\zeta}$, the condition that all $g_{\zeta,\tau} \geq 0$ is equivalent to
\[\sup_{\tau \in \{-1,1\}^{Z_\zeta}} B_\zeta^{-1/2} \left\| \sum_{z \in Z_\zeta} \tau_z \psi_z \right\|_\infty
  = \sup_{\tau \in \{-1,1\}^{Z_\zeta}} \|g_{\zeta,\tau} - \psi_0\|_\infty
  \leq \|\psi_0\|_\infty.\]
% \Bk{Do you mean $\Vert\psi_0\Vert_\infty$ on the rhs?} \CC{Since $\psi_0$ is a constant function, it might be ok as is, with a slight abuse of notation. But $\|\psi_0\|_\infty$ is probably better.}
For the Fourier basis, each $\|\psi_z\|_\infty = 1$,\footnote{This is the only step in the proof that uses any properties specific to the Fourier basis.}
and so
\[\sup_{\tau \in \{-1,1\}^{Z_\zeta}} \left\| \sum_{z \in Z_\zeta} \tau_z \psi_z \right\|_\infty
  \asymp \sup_{\tau \in \{-1,1\}^{Z_\zeta}} \sum_{z \in Z_\zeta} \|\psi_z\|_\infty
  \asymp |Z_\zeta|
  = \zeta^D.\]
Thus, we precisely need $B_\zeta \in \Omega(\zeta^{2D})$, and it is sufficient, for example, that $b_z \in O(\|z\|^{-D/2})$.

\[B_\zeta
  = \sum_{\|z\| \leq \zeta} b_z^{-2}
  \geq \sum_{\|z\| \leq \zeta} \zeta^D
  \asymp \zeta^{2D}.\]

Finally, we show that
\begin{equation}
D_{\text{TV}} \left( \psi_0^n, \frac{1}{2^{|Z_\zeta|}} \sum_{\tau \in \{-1,1\}^{Z_\zeta}} g_{\zeta,\tau}^n \right)
  \leq \frac{1}{2}
\label{ineq:TV_information_bound}
\end{equation}
(where $h^n : \X^n \to [0, \infty)$ denotes the joint likelihood of $n$ IID samples). For any particular $\tau \in \{-1,1\}^{Z_\zeta}$ and $x_1,...,x_n \in \X$, the joint likelihood is
\begin{align*}
g_{\zeta,\tau}^n(x_1,...,x_n)
& = \prod_{i = 1}^n
  \left(
    1 + c_\zeta \sum_{z \in Z_{\zeta}} \tau_z \psi_z(x_i)
  \right) \\
& = 1 + \sum_{\ell = 1}^n \sum_{\substack{i_1,...,i_\ell \in [n] \\ \text{distinct}}} \sum_{j_1,...,j_\ell \in Z_\zeta} \prod_{k = 1}^\ell c_\zeta \tau_{j_k} \phi_{j_k}(x_{i_k}) \\
& = 1 + \sum_{\ell = 1}^n c_\zeta^\ell \sum_{\substack{i_1,...,i_\ell \in [n] \\ \text{distinct}}} \sum_{z_1,...,z_\ell \in Z_\zeta} \prod_{k = 1}^\ell \tau_{z_k} \psi_{z_k}(x_{i_k}).
\end{align*}
Thus, the likelihood of the uniform mixture over $\tau \in \{-1,1\}^{Z_\zeta}$ is
\begin{align*}
& \frac{1}{2^{|Z_\zeta|}} \sum_{\tau \in \{-1,1\}^{Z_\zeta}} g_{\zeta,\tau}^n(x_1,...,x_n) \\
& = 1 + \frac{1}{2^{|Z_\zeta|}} \sum_{\tau \in \{-1,1\}^{Z_\zeta}} \sum_{\ell = 1}^n c_\zeta^\ell \sum_{\substack{i_1,...,i_\ell \in [n] \\ \text{distinct}}} \sum_{z_1,...,z_\ell \in Z_\zeta} \prod_{k = 1}^\ell \tau_{z_k} \psi_{z_k}(x_{i_k}) \\
& = 1
  + \sum_{\ell = 1}^{\lfloor n/2 \rfloor} c_\zeta{2\ell} \sum_{\substack{i_1,...,i_{2\ell} \in [n] \\ \text{distinct}}} \sum_{z_1,...,z_\ell \in Z_\zeta} \prod_{k = 1}^\ell \psi_{z_k}(x_{i_{2k - 1}}) \psi_{z_k}(x_{i_{2k}}),
\end{align*}
where $\lfloor a \rfloor$ denotes the largest integer at most $a \in [0,\infty)$.
This equality holds because, within the sum over $\tau \in \{-1,1\}^{Z_\zeta}$, any term in which any $\tau_z$ appears an odd number of times will cancel. The remaining terms each appear $2^{|Z_\zeta|}$ times. Thus, the total variation distance is
\begin{align}
\notag
& D_{\text{TV}} \left( \psi_0^n, \frac{1}{2^{|Z_\zeta|}} \sum_{\tau \in \{-1,1\}^{Z_\zeta}} g_{\zeta,\tau}^n \right)
  = \frac{1}{2} \left\| \psi_0^n - \frac{1}{2^{|Z_\zeta|}} \sum_{\tau \in \{-1,1\}^{Z_\zeta}} g_{\zeta,\tau}^n \right\|_1 \\
\notag
& = \frac{1}{2} \int_{\X^n} \left| \sum_{\ell = 1}^{\lfloor n/2 \rfloor} c_\zeta^{2\ell} \sum_{\substack{i_1,...,i_{2\ell} \in [n] \\ \text{distinct}}} \sum_{z_1,...,z_\ell \in Z_\zeta} \prod_{k = 1}^\ell \psi_{z_k}(x_{i_{2k - 1}}) \psi_{z_k}(x_{i_{2k}}) \right| \, d(x_1,...,x_n) \\
\label{ineq:big_integral_triangle_ineq}
& \leq \frac{1}{2} \sum_{\ell = 1}^{\lfloor n/2 \rfloor} c_\zeta^{2\ell} \int_{\X^n} \left| \sum_{\substack{i_1,...,i_{2\ell} \in [n] \\ \text{distinct}}} \sum_{z_1,...,z_\ell \in Z_\zeta} \prod_{k = 1}^\ell \psi_{z_k}(x_{i_{2k - 1}}) \psi_{z_k}(x_{2k}) \right| \, d(x_1,...,x_n),
\end{align}
where we used the triangle inequality. By Jensen's inequality (since $\X = [0,1]$),
\begin{align}
\notag
& \int_{\X^n} \left| \sum_{\substack{i_1,...,i_{2\ell} \in [n] \\ \text{distinct}}} \sum_{z_1,...,z_\ell \in Z_\zeta} \prod_{k = 1}^\ell \psi_{z_k}(x_{i_{2k - 1}}) \psi_{z_k}(x_{i_{2k}}) \right| \, d(x_1,...,x_n) \\
\label{ineq:Jensen}
& \leq \sqrt{\int_{\X^n} \left( \sum_{\substack{i_1,...,i_{2\ell} \in [n] \\ \text{distinct}}} \sum_{z_1,...,z_\ell \in Z_\zeta} \prod_{k = 1}^\ell \psi_{z_k}(x_{i_{2k - 1}}) \psi_{z_k}(x_{i_{2k}}) \right)^2 \, d(x_1,...,x_n)}.
\end{align}
Since $\{\psi_z\}_{z \in \Z}$ is an orthogonal system in $\L^2(\X)$, we can pull the summations outside the square, so
\begin{align*}
& \int_{\X^n} \left( \sum_{\substack{i_1,...,i_{2\ell} \in [n] \\ \text{distinct}}} \sum_{z_1,...,z_\ell \in Z_\zeta} \prod_{k = 1}^\ell \psi_{z_k}(x_{i_{2k - 1}}) \psi_{z_k}(x_{i_{2k}}) \right)^2 \, d(x_1,...,x_n) \\
& = \sum_{\substack{i_1,...,i_{2\ell} \in [n] \\ \text{distinct}}} \sum_{z_1,...,z_\ell \in Z_\zeta} \int_{\X^n} \left( \prod_{k = 1}^\ell \psi_{z_k}(x_{i_{2k - 1}}) \psi_{z_k}(x_{i_{2k}}) \right)^2 \, d(x_1,...,x_n) \\
& = \sum_{\substack{i_1,...,i_{2\ell} \in [n] \\ \text{distinct}}} \sum_{z_1,...,z_\ell \in Z_\zeta} 1
  = \binom{n}{2\ell} \zeta^{D\ell}
  \leq \frac{n^{2\ell} \zeta^{D\ell}}{(\ell!)^2},
\end{align*}
since
\[\binom{n}{2\ell}
  = \frac{n!}{(2\ell)!(n - 2k)!}
  \leq \frac{n^{2\ell}}{(2\ell)!}
  \leq \frac{n^{2\ell}}{(\ell!)^2}.\]
Combining this with inequalities \eqref{ineq:big_integral_triangle_ineq} and \eqref{ineq:Jensen} gives
\begin{equation}
D_{\text{TV}} \left( \psi_0^n, \frac{1}{2^{|Z_\zeta|}} \sum_{\tau \in \{-1,1\}^{Z_\zeta}} g_{\zeta,\tau}^n \right)
  \leq \frac{1}{2} \sum_{\ell = 1}^{\lfloor n/2 \rfloor} \frac{\left( n c_\zeta^2 \zeta^{D/2} \right)^\ell}{\ell!}
  \leq \exp \left( n c_\zeta^2 \zeta^{D/2} \right) - 1,
\label{ineq:exponential_bound}
\end{equation}
where we used the fact that the exponential function is greater than any of its Taylor approximations on $[0, \infty)$. The last expression in inequality~\eqref{ineq:exponential_bound} vanishes if $n c_\zeta^2 \zeta^{D/2} \to 0$. Recalling now that we set $c_\zeta = B_\zeta^{-1/2}$, for some constant $C > 0$, the desired bound \eqref{ineq:TV_information_bound} holds by choosing $\zeta$ satisfying
\[\frac{\zeta^D}{B_\zeta^2}
  = \zeta^D c_\zeta^4
  \leq C n^{-2}.\]

{\bf The unsmooth case ($B_\zeta \in o \left( \zeta^{2D} \right)$):} Finally, we consider the `highly unsmooth' case, when $B_\zeta \in o(\zeta^{2D})$ 
% \Bk{u mean $\zeta$ instead of $z$?} \CC{Yes, fixed.}.
In this case, we must modify the above proof to ensure that the $g_{\zeta,\tau}$ functions are all non-negative. In the Fourier case, we again wish to ensure
\[c_\zeta \zeta^D
  = c_\zeta |Z_\zeta|
  \asymp c_\zeta \sup_{\tau \in \{-1,1\}^{Z_\zeta}} \left\| \sum_{z \in Z_\zeta} \tau_z \psi_z \right\|_\infty \leq 1,\]
but this is no longer guaranteed by setting $c_\zeta = B_\zeta^{-1/2}$; instead, we use the smaller value $c_\zeta = \zeta^{-D}$. Clearly, we still have $\|g_{\zeta,\tau}\|_b^2 \leq 1$. Now, however, we have a smaller estimation error
\begin{equation}
\inf_{\tau \in \{-1,1\}^{Z_\zeta}} \|g_{\zeta,\tau}\|_a^2 - \|\psi_0\|_a^2
  = c_\zeta^2 \sum_{z \in Z_\zeta} a_z^{-2}
  = \frac{A_\zeta}{\zeta^{2D}}.
\label{eq:smaller_gap}
\end{equation}
Also, the information bound~\eqref{ineq:exponential_bound} now vanishes when $n \zeta^{-3D/2} = n c_\zeta^2 \zeta^{D/2} \to 0$, so that, for some constant $C > 0$, the desired bound \eqref{ineq:TV_information_bound} holds by choosing $\zeta$ satisfying
\[\zeta \leq C n^{2/(3D)}.\]
Plugging this into equation~\eqref{eq:smaller_gap} gives
\[\inf_{\tau \in \{-1,1\}^{Z_\zeta}} \|g_{\zeta,\tau}\|_a^2 - \|\psi_0\|_a^2
  \asymp \frac{A_{n^{2/(3D)}}}{n^{4/3}}.\]
Finally, by Le Cam's lemma, this implies the minimax rate
\[\inf_{\hat S} \sup_{p,q \in \H_b} \E \left[ \left( \hat S - \langle p, q \rangle \right)^2 \right]
  \geq \left( \frac{A_{n^{2/(3D)}}}{n^{4/3}} \right)^2.\]
% Note that, $A_\zeta \asymp O \left( \zeta^{2s + D} \right)$, this implies
% \[\inf_{\hat S} \sup_{p,q \in \H_b} \E \left[ \left( \hat S - \langle p, q \rangle \right)^2 \right]
%   \geq \frac{n^{\frac{8s + 4D}{3D}}}{n^{8D/(3D)}}
%   = n^{\frac{8s - 4D}{3D}}.\]
% When is $\frac{8s - 4D}{3D} \geq -1$? If $D \leq 8s$, this is slower than $n\inv$.
\end{proof}

\section*{Acknowledgements}
SS is supported by a National Science Foundation Graduate Research Fellowship
under Grant No. DGE-1252522. BKS is supported by NSF-DMS-1713011.

\bibliographystyle{abbrvnat}
\bibliography{biblio}

\end{document}